\documentclass[notitlepage,oneside]{report}

\usepackage{latexsym,amscd,amssymb,amstext,epic,curves}

\begin{document}

\newcommand{\mcross}{\begin{picture}(10,10)
                        \put(0,0){\line(1,1){4.3}}
                        \put(10,10){\line(-1,-1){4.3}}
                        \put(0,10){\line(1,-1){10}}
                        \put(0,10){\vector(-1,1){0}}
                        \put(10,10){\vector(1,1){0}}
                        \end{picture}   }

\newcommand{\pcross}{\begin{picture}(10,10)
                        \put(0,0){\vector(1,1){10}}
                        \put(10,0){\line(-1,1){4.3}}
                        \put(0,10){\line(1,-1){4.3}}
                        \put(0,10){\vector(-1,1){0}}
                        \end{picture}   }

\newcommand{\ocross}{\begin{picture}(10,10)
                       \put(0,5){\oval(7,10)[r]}
                        \put(10,5){\oval(7,10)[l]}
                        \put(10,10){\vector(1,0){0}}
                        \put(0,10){\vector(-1,0){0}}
                        \end{picture}   }

\vspace{-2cm}

\title{Presentations of homotopy skein modules of oriented $3$-manifolds}

\makeatletter
\date{Uwe Kaiser \footnote{Partially supported by DFG 653/2000}\\
        \ \\
        \begin{small}Universit\"{a}t Siegen \\
        57068 Siegen, Germany \\
        kaiser@mathematik.uni-siegen.de \end{small} \\
         \   \\
         July 6, 2000}
\makeatother

\maketitle

\begin{small}
\begin{quote}
\textsc{ABSTRACT.} \ A new method to derive presentations of
skein modules is developed. For the case of homotopy skein
modules it will be shown how the topology of a $3$-manifold is
reflected in the structure of the module. The freeness problem
for $q$-homotopy skein modules is solved, and a natural skein
module related to linking numbers is computed.
\end{quote}
\end{small}

\vspace{0.5cm}

\centerline{\textsc{INTRODUCTION}}

\vspace{0.3cm}

Let $R$ be a commutative ring. A skein module of an oriented
$3$-manifold $M$ is defined from the free $R$-module on a set of
isotopy classes of links (possibly decorated by orientation or
framing) in $M$ by dividing out a submodule generated by local
skein relations [P-1]. The skein modules of $3$-manifolds
$F\times I$ for $F$ a compact oriented surface play a central
role in the study of quantum invariants, topological quantum field
theory and finite type invariants of links in $3$-manifolds. For
important skein relations it is known that the modules of
$F\times I$ can be derived through quantization of algebras of
homotopy classes of loops on $F$ [B-F-K, H-P, P-3, P-4, T].
Moreover several recent results have shown that skein modules
of $3$-manifolds should be interpreted as deformations of algebras
emerging from the theory of linear representations of the fundamental group
[B, P-S]. A theory of skein modules for arbitrary $3$-manifolds
\textit{beyond} $F\times I$ seems to be crucial in understanding
the quantum topology of $3$-dimensional manifolds. But it is
difficult to derive global results about skein modules, and
complete structure results are known only for simple skein
relations [P-1, P-5]. A workable description of the structure of a
class of skein modules always has to face two different problems:

\vspace{0.3cm}

\noindent \textbf{Problem 1:} Find a suitable set of generators for the skein module.

\vspace{0.2cm}

\noindent \textbf{Problem 2:} Describe the relations with respect
to this generating set in terms of the topology of $M$.

\vspace{0.3cm}

Problem 1 is related with subtle finiteness questions concerning
topology and combinatorics of links in $3$-manifolds [Ka-L]. With
respect to \textit{minimal} generating sets, Problem 2 asks how
the topology obstructs the existence of certain skein invariants.
It is known that skein invariants are related to finite type
invariants, which have been studied intensively.  Thus it is
natural to study problem 2 modulo problem 1. This has been
carried out in various approaches in terms of (completed) isotopy
skein modules [B-F-K, Ka-2].

\vspace{0.2cm}

We will show that a direct approach to \textit{global} results
concerning skein modules is possible along the ideas of [Ka-L],
[K] and [Ki-Li]. Our main philosophy is: The obstructions to the
existence of finite type invariants are relations in a
presentation of skein modules. For homotopy skein modules problem
1 has a natural and easy answer. So it is possible to focus on
problem 2 from scratch (without completion). We will give a quite
complete general answer to problem 2 for this case, though often
it is still difficult to solve explicit problems by our method. A
discussion of presentations for the (completed) universal
Jones-Conway skein module as defined in [T, P-4] and for Kauffman
skein modules [P-1, H-P] will be left to future work.

\vspace{0.3cm}

We like to thank Brandeis University and George Washington
University for their hospitality, and Jerry Levine,
Jozef Przytycki, Danny Ruberman and Manfred Stelzer
for many helpful discussions and suggestions.

\vspace{0.5cm}

\centerline{\textsc{1.\ STATEMENT OF THE MAIN RESULTS}}

\vspace{0.3cm}

Throughout let $R=\mathbb{Z}[q^{\pm 1},z]$ and let $M$ be a compact
oriented $3$-manifold. Let $\mathfrak{H}(M)$ denote the set of
link homotopy classes of oriented links in $M$, including the
empty link $\emptyset $.  Two oriented links are
\textit{link homotopic} if they are homotopic through a deformation, which
keeps different components disjoint.

\vspace{0.2cm}

Following J.\ Przytycki [P-2]
the $q$-homotopy skein module $\mathcal{H}(M)$ is the quotient of
$R\mathfrak{H}(M)$ by the submodule, which is generated by all skein elements
$q^{-1}K_+-qK_--zK_0$
for all crossings of different components of links (i.\ e.\
\textit{mixed crossings}).

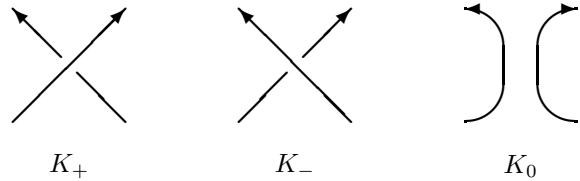
\begin{figure}[htbp]
        \setlength{\unitlength}{1.5mm}
        \thicklines
        \centering
        \begin{picture}(50,10)

                \put( 0, 0){\pcross}
                \put(20, 0){\mcross}
                \put(40, 0){\ocross}

                \put( 5,-4){\makebox(0,0){$K_+$}}
                \put(25,-4){\makebox(0,0){$K_-$}}
                \put(45,-4){\makebox(0,0){$K_0$}}

        \end{picture}
\end{figure}

\vspace{1cm}

The modules $\mathcal{H}(M)$ are universal with respect to the
Homfly homotopy skein relation [P-3]. There are two obvious ways
to simplify the relation. Let $\mathcal{C}(M)$ resp.\
$\mathcal{L}(M)$ denote the $\mathbb{Z}[z]$- resp.\
$\mathbb{Z}[q^{\pm 1}]$-modules resulting from $\mathcal{H}(M)$ by
the ring homomorphisms $q\mapsto 1$ resp.\ $z\mapsto 0$. It
follows from Przytycki`s Universal Coefficient Theorem [P-3] that
$\mathcal{C}(M)\cong \mathcal{H}(M)\otimes_R \mathbb{Z}[z]$ resp.\
$\mathcal{L}(M)\cong \mathcal{H}(M)\otimes_R \mathbb{Z}[q^{\pm 1}]$ with
the $R$-module structures on $\mathbb{Z}[z]$ and
$\mathbb{Z}[q^{\pm 1}]$ given by the ring homomorphisms above. In this way results about
$\mathcal{H}(M)$ also give rise to results about the \textit{derived}
skein modules.

\vspace{0.2cm}

For $M$ a connected $3$-manifold, $\mathcal{H}(M)\cong R\oplus
\tilde{\mathcal{H}}(M)$, where $\tilde{\mathcal{H}}(M)$ is the
module defined from $\mathfrak{H}(M)\setminus \emptyset$. In
general, $\mathcal{H}(M)\cong R\oplus \oplus_i
\tilde{\mathcal{H}}(M_i)$, where the sum is over the set of path
components $M_i$ of $M$. Throughout the following we assume that
$M$ is connected and we fix a basepoint $*$. All homotopy groups
will be defined with respect to $*$. Let $\hat{\pi }(M)$
denote the set of conjugacy classes of elements of $\pi_1(M)$. Let
$\mathfrak{b}(M)$ denote the set of unordered sequences of
elements of $\hat{\pi }(M)$ of length $\geq 0$. The surjective
\textit{wrapping invariant}
$$\omega :\mathfrak{H}(M)\rightarrow \mathfrak{b}(M)$$
is defined by assigning to each oriented link in $M$ the
unordered sequence of free homotopy classes of its components.

\vspace{0.2cm}

A section $\sigma $ of $\omega $ defines a set of \textit{standard
links} $K_{\alpha }:=\sigma (\alpha )$, and induces the
epimorphism of $R$-modules
$$\sigma : SR\hat{\pi }(M)\cong R\mathfrak{b}(M)\twoheadrightarrow \mathcal{H}(M).$$
Here, for each module $A$ we let $SA$ denote the symmetric
algebra on $A$. (We let $\sigma $ also denote the corresponding
$\mathbb{Z}[z]$-epimorphism onto $\mathcal{C}(M)$.)
The surjectivity of the homomorphism $\sigma $ is easily proved
by induction on the number of components of links. It follows from the
simple observation that smoothing a mixed crossing reduces this
number.

\vspace{0.2cm}

A set of generators $\mathfrak{G}$ of a module $H$ is
\textit{minimal}, if for each $g\in \mathfrak{G}$ the set
$\mathfrak{G}\setminus g$ does not generate $H$. The image of
$\sigma $ is a \textit{minimal} set of generators of
$\mathcal{H}(M)$ (and similarly for $\mathcal{C}(M)$ and
$\mathcal{L}(M)$). This can be proved by applying the ring
homomorphism $q\mapsto 1, z\mapsto 0$. It maps
$\mathcal{H}(M)$ onto the free abelian group generated by the set
of oriented links in $M$ with relations $K_+=K_-$ for \textit{all}
crossings. But this free abelian group is isomorphic to
$S\mathbb{Z}\hat{\pi }(M)$.

\vspace{0.2cm}

Let $\hat{M}$ be the $3$-manifold defined from $M$ by capping off
all $2$-spheres in the boundary by standard $3$-balls. It is
well-known [P-1] that the inclusion $M\rightarrow \hat{M}$ induces
isomorphisms of skein modules. A $3$-manifold $M$ is
\textit{atoroidal} [J-S] if each map of a torus $S^1\times
S^1\rightarrow M$, which induces an injective homomorphism of
fundamental groups, can be homotoped into the boundary of $M$.
Let $b_1(X)$ denote the first Betti number of a complex $X$.

\vspace{0.2cm}

\noindent \textbf{Theorem 1.1.}

\noindent a) $\mathcal{H}(M)$ \textit{is free if and only if
$\pi_1(M)$ is abelian and $2b_1(M)=b_1(\partial M)$. Then
$\sigma $ is an isomorphism. Otherwise $\mathcal{H}(M)$ has
torsion.}

\noindent b) \textit{If $\pi_2(\hat{M})=0$ and $M$ (or
$\hat{M}$) is atoroidal then $\mathcal{C}(M)$ is free (and
$\sigma $ is an isomorphism).}

\vspace{0.3cm}

Theorem 1.1 generalizes the results of J.\ Przytycki [P-3] and J.\
Hoste and J.\ Przytycki [H-P] for the case $M=F\times I$ and $F$ a
compact oriented surface. Our proofs do not rely on their methods
or results but are derived by conceptually new and independent
ideas. In particular we prove freeness of $\mathcal{C}(M)$ for all
$3$-manifolds with complete hyperbolic structures. (It is known
that the modules $\mathcal{C}(M)$ are free if and only if they
are torsion-free [K-1].) It also follows that $q$-homotopy skein
modules of Lens spaces are free as conjectured by J.\ Przytycki in
[P3].

\vspace{0.2cm}

The oriented intersection number between oriented closed surfaces and
loops defines the intersection pairing (see [K-2] for details):
$$\iota :H_2(M)\otimes H_1(M)\rightarrow \mathbb{Z}.$$

Let $LM$ denote the free loop space of $M$ and let
$f:S^1\rightarrow M$ be a basepoint. Elements of $\pi_1(LM,f)$
can represented by maps $S^1\times S^1\rightarrow M$. The map $\iota $ composes with natural
maps $\pi_1(LM,f)\rightarrow H_2(M)$ and
$\mathfrak{b}(M)\rightarrow H_1(M)$ (adding the homology classes of the elements of
a sequence)
to define
$$\iota_f: \pi_1(LM,f)\times \mathfrak{b}(M)\longrightarrow \mathbb{Z}.$$

Let $f$ represent some $a\in \hat{\pi }(M)$, which appears in $\alpha $.
Let $\alpha \setminus a$ be defined by omitting some element $a$
from the sequence. Then let $\iota (\alpha,a)$ be the positive
generator of the subgroup $\iota_f( \pi_1(LM,f), \alpha \setminus
a)\subset \mathbb{Z}$. Let $\lambda (\alpha )$ be the greatest
common divisor of all $\iota (\alpha,a)$ for all $a$ in $\alpha $.

\vspace{0.2cm}

\noindent \textbf{Theorem 1.2.} \ \textit{The module
$\mathcal{L}(M)$ is isomorphic to  $$\bigoplus \limits_{\alpha
\in \mathfrak{b}(M)}\mathbb{Z}[q^{\pm 1}]/(q^{2\lambda (\alpha
)}-1)\mathbb{Z}[q^{\pm 1}].$$}

The structure of $\mathcal{L}(M)$ can be interpreted as a
statement about \textit{linking numbers}. In fact, by 1.2 we can
define \textit{linking numbers}(relative to the choice of
standard links) $\ell_{\alpha }: \omega^{-1}(\alpha )\rightarrow
\mathbb{Z}/\lambda (\alpha ) \mathbb{Z}$ satisfying the property:
$\ell_{\alpha }$ changes by $\pm 1$ through $\pm 1$-crossings of
different components of a link. For given $\alpha $, the linking
number is defined in $\mathbb{Z}$, if and only if intersection
numbers of singular tori, defined by self-homotopies of a
component of $\alpha $, with the remaining components of $\alpha
$, vanish. It follows that linking numbers are globally defined,
or equivalently, $\mathcal{L}(M)$ is free if and only if singular
tori in $M$ are homologous into the boundary of $M$ (see also [K-2]).

\vspace{0.3cm}

The results above will be deduced from a general presentation for
homotopy skein modules, which we describe now. A set of standard
links $\{K_{\alpha }\}$ is \textit{geometric} if the following
conditions hold for each $\alpha $:

\begin{enumerate}
\item $K_{\alpha }$ is an oriented embedded link in $M$.
\item Each sublink of $K_{\alpha }$ is a standard link.
\item If $a^{\pm 1}\in \hat{\pi }(M)$ appears $k$ times in $\alpha $ then $K_{\alpha }$ contains
$k$ parallel copies of a knot in $M$  with the suitable
orientations.(Here $a^{-1}$ is the inverse conjuagacy class, and
parallel refers to some framing).
\item The sublink of null-homotopic components of $K_{\alpha }$ is a link with unknotted and
unlinked components contained in a $3$-ball in $M$.
\end{enumerate}

It is not hard to prove that each $3$-manifold admits a geometric
set of standard links. In fact, in a handle-body such a set can
easily be constructed by using monotonicity with respect to some
$I$-structure. For general $M$ choose a Heegard decomposition.
Then consider a suitable subset of a set of geometric standard
links for the skein module of one of the handle-bodies $F$ of the decomposition. The
subset can be defined from a section of the surjection $\pi
(F)\rightarrow \pi(M)$.

\vspace{0.2cm}

Given a geometric set of standard links we choose orderings and
basings of its components. The basings are paths, which join
$*\in M$ with points on the components. Because of condition 2 in
the definition above we can chose the basings compatible up to
isotopy with sublinks. Thus in particular we have chosen a
section of the projection map $\pi_1(M)\rightarrow \hat{\pi }(M)$ and
an ordering of the sequences in $\mathfrak{b}(M)$, which will
apply to all following constructions.

\vspace{0.2cm}

\noindent {$\Delta$-\textbf{Construction.} \ Let $c =(a,b)$ be a
pair of elements of $\pi_1(M)$, which appears in $\alpha $ (in
some order). Let $K_{\alpha }^c$ be the oriented link, which
results from $K_{\alpha}$ by connecting two components with
homotopy classes $a,b$ by a band, which follows the given basings.
We can change the homotopy class of the band by a loop in $*$ representing
some $g\in \pi_1(M)$. This defines a new banded link
$K_{\alpha,g}^c$. The difference $K_{\alpha,g}^c-K_{\alpha }^c$
can be expanded (not uniquely) in terms of standard links with
less components than $K_{\alpha }$. Let $\tilde{\Delta }(\alpha
,c,g)\in SR\hat{\pi }(M)$ be some element, which results from such an
expansion. Let
$\Delta(\alpha,c,g):=z(q^{-1}-q)\tilde{\Delta}(\alpha,c,g)$.

\vspace{0.3cm}

\noindent $\Theta$-\textbf{Construction.} \ Let $a\in \pi_1(M)$
appear in $\alpha $ and let $f: S^1\rightarrow M$ be a
parametrized component of $K_{\alpha }$ with homotopy class $a$.
Each element $h\in \pi_1(LM,f)$ defines a self-homotopy $f_t$ of
$f$ in $M$ (which defines a singular torus in $M$). We can assume
that $f_t$ intersects $K_{\alpha ,a}$ for only a finite number of
$t$ transversely and in a single point for each $t$. Here
$K_{\alpha ,a}$ is the link, which results by forgetting a
component with homotopy class $a$. Each intersection of $f_t$
with $K_{\alpha ,a}$ comes naturally with a sign and a smoothing
thus defining a linear combination of links. We can expand each
term in this linear combination in terms of standard links, again
not uniquely of course, to define some element
$\Theta(\alpha,a,h)\in SR\hat{\pi }(M)$.

\vspace{0.3cm}

\noindent \textbf{Theorem 1.3.} \textit{Let $\sigma $ define a
geometric set of standard links. Then the kernel of $\sigma $ is
generated by the elements $\Delta (\alpha ,c,g)$ and $\Theta
(\alpha ,a,g)$, which result from the
$\Delta$/$\Theta$-construction.}

\vspace{0.3cm}

By 1.3, the relations of the homotopy skein module are
parametrized by copies of the fundamental group of $M$ and of its
free loop spaces (singular tori in $M$). This reduces problem 2
of the introduction to questions about curves and singular tori
in $3$-manifolds. Though the homotopy classification of singular
tori in $3$-manifolds is well-understood (see e.\ g.\ [J-S]),
explicit problems about skein modules can still be difficult to
handle. In [K-4] we discuss how torsion problems of homotopy
skein modules are interpreted in terms of
\textit{non-commutative} intersection maps.

\vspace{0.2cm}

Theorem 1.3 will be deduced from a skein theory of paths in mapping spaces with restricted singularities, which is the central point of this paper.

\vspace{0.3cm}

Here is the plan for the rest of the paper:

\vspace{0.2cm}

In section 2 we state two refined versions of theorem 1.3:
Presentations for homotopy skein modules defined from restricted
sets of links, and presentations of $\mathcal{H}(M)$ as module
over the homotopy skein algebra of $S^3$. The first presentations
are interesting from the viewpoint of a study of link homotopy
invariants in $3$-manifolds. In section 3 we state the necessary
background about mapping spaces related with the Lin approach
and show how the ideas apply to skein theory. In
section 4 we describe our basic construction of a skein theory of
paths. This is applied in section 5 to prove theorem 1.3. In
section 6 we develop some homotopy theory of free loop spaces. In
section 7 we show that relations given by the
$\Delta$/$\Theta$-construction vanish for many situations. This
allows to conclude theorems 1.1 and 1.2 in section 8.

\vspace{0.5cm}

\centerline{\textsc{2.\ THE REFINED PRESENTATION RESULTS}}

\vspace{0.3cm}

Following [K-1] we introduce homotopy skein modules over restricted
sets of links. Let $K=K_+$ be a link with $\omega (K_{\pm})=
\alpha $. Then we describe the wrapping invariant of a smoothing
$K_0$ in the following way: Choose an ordering and basing of the
components of $K$. This induces an ordering of $\omega (K)
=\langle \alpha_1,\alpha_2, \ldots ,\alpha_r\rangle$ to a sequence
$(a_1,a_2,\ldots ,a_r)$ of elements of $\pi_1(M)$. Consider the
ball in $M$, where two components of $K$ are smoothed. Let us
assume that these are the first two components and the
overcrossing arc comes from the first component. Choose a path
$\gamma $, which joins $*\in M$ to a point in the boundary of the
$3$-ball. The homotopy class of the smoothed component (with new
basepoint on the right arc of the smoothed link in the ball) is
given by $g_1a_1g_1^{-1}g_2a_2g_2^{-1}$ for suitable elements
$g_i\in \pi_1(M)$ and $i=1,2$. E.\ g.\ the element $g_1$ is
defined by juxtaposing $\gamma $ with a path in the smoothing
ball to the crossing point, then with the segment on $K_1$ from
the crossing point on $K_1$ to the basepoint $*_1$ on $K_1$
(opposite to the orientation on $K_1$), and finally with the
basing arc for $K_1$ from $*_1$ back to $*$. The corresponding
wrapping invariant
$$\langle (a_1ga_2g^{-1})^{\circ},\alpha_3,\ldots ,\alpha_r)$$
with $g:=g_1^{-1}g_2$, is called a (first order)
\textit{descendant} of $\alpha $. Here $a^{\circ}$ denotes the
conjugacy class of $a\in \pi_1(M)$. By iteration this defines the
\textit{set of descendants} of $\alpha $.

\vspace{0.2cm}

A subset $C\subset \mathfrak{b}(M)$ is called \textit{skein
closed} if it contains the descendants of its elements. The
\textit{skein closure} $\bar{C}$ of $C$ is the intersection of all
skein closed sets containing $C$. Let
$\mathfrak{H}_C(M):=\omega^{-1}(C)$. We define the $q$-homotopy
skein module $\mathcal{H}_C(M)$ to be the quotient of
$R\mathfrak{H}_C(M)$ by the submodule generated by all skein
relations $q^{-1}K_+-qK_--zK_0$ for all $K_+\in \omega^{-1}(\bar{
C})$. It follows from the definitions that $\sigma $ restricts to
epimorphisms $R\bar{C}\rightarrow \mathcal{H}_C(M)$.

\vspace{0.2cm}

We need further notation. Let $\mathfrak{b}_r(M)$ denote the set
of elements $\alpha \in \mathfrak{b}(M)$ of fixed length
$r=:|\alpha |$. The unique sequence of length $0$ is denoted
$\langle \rangle $. For $\alpha \in \mathfrak{b}_r(M)$ and a
subsequence $\alpha ' \in \mathfrak{b}_s(M)$ let $\alpha
\setminus \alpha ' \in \mathfrak{b}_{r-s}(M)$ be the sequence
defined by eliminating the elements of $\alpha ' $ from $\alpha $.
For $C$ skein closed let $C^{\ast}$ be the set of those
sequences, which result from adjoining arbitrary trivial
conjugacy classes to sequences of $C$. The sets $C^{\ast}$ are
obviously skein closed. Let $\mathfrak{t}(M):=\langle
\rangle^{\ast}$ be the set of all sequences with only trivial
conjucagy classes.

\vspace{0.2cm}

For $g\in \pi_1(M)$ let $\langle g\rangle$ denote the cyclic
subgroup generated by $g$. Given $a,b\in \pi_1(M)$, let $D(a,b)$
denote the set of non-trivial double left $\langle a \rangle $
and right $\langle b \rangle$ cosets of elements of $\pi_1(M)$.
For $a\in \pi_1(M)$ let $T(a)$ denote a set of generators of
$\pi_1(LM,f)$ for $f$ a parametrized component of standard links
representing $a$. For $\alpha \in C$, and $(a,b)$ resp.\ $a$
contained in $\alpha $, the $\Delta$- resp.\
$\Theta$-constructions define maps $D(a,b)\rightarrow SR\bar{C}$
resp.\ $T(a)\rightarrow SR\bar{C}$. (It will be shown in section
5 that the $\Delta$-construction is defined on cosets as above.)
For each $\alpha \in C$ we define the \textit{structure set}:
$$S(\alpha ):= \bigcup_{a}T(a)\cup \bigcup_{(a,b)}D(a,b),$$
with the union over all distinct $a$ in $\alpha $ resp.\ distinct
pairs $(a,b)$ in $\alpha $ with both $a,b$ non-trivial, such that
$\alpha \setminus  a^{\circ} $ resp.\ resp.\ $\alpha \setminus
\langle a^{\circ},b^{\circ}\rangle $ are not contained in
$\mathfrak{t}(M)$. The maps so defined combine to define
$$\chi_{\alpha }: \mathcal{S}(\alpha )\longrightarrow SR\bar{C}$$
by $\chi_{\alpha }(g)=\Delta (\alpha ,(a,b),g)$ for $g\in D(a,b)$,
and $\chi_{\alpha }(h)=\Theta (\alpha ,a,h)$ for $h\in T(a)$. So
we have defined the homomorphism
$$\chi =\sum \chi_{\alpha}: \bigoplus \limits_{\alpha \in \bar{C}}R\mathcal{S}(\alpha )
\longrightarrow R\bar{C}.$$ The definition involves choices of
expansions of links and is not defined in a unique way.

\vspace{0.2cm}

\noindent \textbf{Theorem 2.1.} \textit{Let $C$ be a subset of
$\mathfrak{b}(M)$ and let $\sigma $ define a geometric set of
standard links. Then the sequence
$$\CD \bigoplus \limits_{\alpha \in \bar{C}}R\mathcal{S}(\alpha )@>{\chi }>>
R\bar{C} @>{\sigma }>> \mathcal{H}_C(M)\rightarrow 0
\endCD$$
is exact for all choices of $\chi $ determined by the $\Delta$/$\Theta$-construction.
}

\vspace{0.3cm}

For $C=\mathfrak{b}(M)$ theorem 2.1 reduces to 1.3.

\vspace{0.2cm}

\noindent \textbf{Example.} \ For each $3$-manifold $M$ the module
$\mathcal{H}_{\mathfrak{t}(M)}(M)$ is free. On the other hand it
is not known whether in general the image of
$R\mathfrak{H}_{\mathfrak{t}(M)}(M)$ in $\mathcal{H}(M)$, which
is also the natural image of $\mathcal{H}(D^3)$ in
$\mathcal{H}(M)$, is free.

\vspace{0.2cm}

\noindent \textbf{Remark 2.2.} \ The above presentations are not
always natural with respect to oriented embeddings of
$3$-manifolds $j: M\hookrightarrow N$. Let $j_*:
\mathcal{H}(M)\rightarrow \mathcal{H}(N)$ be the homomorphism
induced by $j$. Let $\sigma_M$ and $\sigma_N$ be any choices of
standard links. Then there is a homomorphism
$c(\sigma_M,\sigma_N)$, which makes the following diagram commute:
$$
\CD
SR\hat{\pi }(M)@>{\sigma_M}>>\mathcal{H}(M)\\
@V{c(\sigma_M,\sigma_N)}VV  @VV{j_*}V \\
SR\hat{\pi }(N)@>\sigma_N>>\mathcal{H}(N)
\endCD
$$
The homomorphism $c(\sigma_M,\sigma_N)$ expands the standard
links given by $\sigma_M$ in terms of standard links from
$\sigma_N$. But in general this  homomorphism is \textit{not}
induced from the corresponding induced maps on fundamental groups
and sets of conjugacy classes. On the other hand, let $j_{\sharp}:
\hat{\pi }(M)\rightarrow \hat{\pi }(N)$ be injective. The $j$ induces
injective maps $j_*: \mathfrak{b}(M) \rightarrow \mathfrak{b}(N)$
and, for $\alpha \in \mathfrak{b}(M)$, injective maps $j_*:
\mathcal{S}(\alpha )\rightarrow \mathcal{S}(j_*(\alpha ))$.
Finally we have induced homomorphisms
$$j_*: \bigoplus \limits_{\alpha \in C}R \mathcal{S}(\alpha )\rightarrow \bigoplus
\limits_{\beta \in (j(C))}R\mathcal{S}(\beta )$$ for $C\subset
\mathfrak{b}(M)$ skein closed. In this case it is easy to choose
\textit{geometric} standard links such that
$c(\sigma_M,\sigma_N)=j_{\sharp }$. There is the commutative
diagram for all sets $C\subset \mathfrak{b}(M)$ and
$j_*(C)\subset D\subset \mathfrak{b}(N)$:
$$\CD
\bigoplus \limits_{\alpha \in (C)}\mathcal{S}(\alpha )@>{\chi_M
}>>SR\bar{C}@>{\sigma_M}>> \mathcal{H}_C(M)@>>>0
\\ @V{j_*}VV @V{j_*}VV @V{j_*}VV \\
\bigoplus \limits_{\beta \in (D)}\mathcal{S}(\beta )@>{\chi_N
}>>SR(D)@>{\sigma_N}>>\mathcal{H}_D(N)@>>>0
\endCD$$
This applies in particular to $C=\mathfrak{b}(M)$ and
$D=\mathfrak{b}(N)$, but also to $M=N$ and $C\subset D$.

\vspace{0.3cm}

It is immediate from 2.1 that $\mathcal{H}(M)$ is a free module
for each homotopy $3$-sphere (by 6.1.\ or [L], for $M$ a homotopy
sphere, $\pi_1(LM,f)$ is the trivial group for each map
$f:S^1\rightarrow M$). For $S^3$ the result has been proven by
J.\ Przytycki [P-3]. Moreover, $\sigma $ is an isomorphism and
$\mathcal{H}(S^3)\cong \mathcal{H}(D^3)$ is isomophic to the
polynomial algebra:
$$\mathcal{R}:=\mathbb{Z}[q^{\pm 1},z,u],$$ the multiplication
defined by stacking links.

\vspace{0.2cm}

We want to describe $\mathcal{H}(M)$ as module over $\mathcal{H}(D^3)$.

\vspace{0.2cm}

\noindent \textbf{Lemma.} \ \textit{The inclusion of an oriented $3$-ball in $M$
induces the canonical map
$$\mathfrak{H}(D^3)\times \mathfrak{H}(M)\rightarrow \mathfrak{H}(M),$$
which defines the strucure of a $\mathcal{R}$-module on $\mathcal{H}(M)$,
i.e. a homomorphism of $R$-modules:
$$\mathcal{H}(D^3)\otimes \mathcal{H}(M)\rightarrow \mathcal{H}(M).$$}

\vspace{-0.3cm}

\noindent \textsl{Proof.} \ We can replace $M$ by $M\setminus int(D^3)$.
Then define the map of link homotopy sets by forming  the union of $D^3$ and $M\setminus
int(D^3)$.  The homomorphism of $R$-modules is deduced in the obvious way. Since any two
$3$-balls in $M$ are isotopic the choice of $3$-ball does not change the structure of module.
$\square$

\vspace{0.2cm}

\noindent \textbf{Remark.} \ The $\mathcal{R}$-module $\mathcal{H}(M)$ is
$\mathcal{R}$-isomorphic to the following skein module:
The generating set is given by the set of all link homotopy classes of oriented links in $M$ and
the relations are given by
\begin{enumerate}
\item[i)]
$q^{-1}K_+-qK_--zK_0$ for all skein triples, and
\item[ii)] $tK-(K\coprod U)$, where $U$ is the unknot contained in a $3$-ball in $M$
separated from $K$.
\end{enumerate}

\vspace{0.2cm}

Recall that the
sublink of trivial components of each link $K_{\alpha}$ is a
standard unlink contained in a $3$-ball in $M$. Thus $K_{\alpha }$
is the product of this unlink and some standard link $K_{\alpha
'}$ with only non-trivial free homotopy classes in $\alpha '$.  It
follows that the $\mathcal{R}$-module $\mathcal{H}(M)$ hs
generated by standard links indexed by $\hat{\mathfrak{b}}(M)=\mathfrak{b}(M)\setminus \mathfrak{t}(M)$.

\vspace{0.2cm}

Now we define the \textit{structure sets}
$\mathcal{S}_{\mathcal{R}}(\alpha )$ for the presentations over
$\mathcal{R}$. Let $\mathfrak{S}$ be a set
of generators of $\pi_2(M)$. For $\alpha \in \hat{\mathfrak{b}}(M)$ let
$$\mathcal{S}_{\mathcal{R}}(\alpha):=\mathfrak{S} \cup \bigcup_a T(a)\cup
\bigcup_{(a,b)}D(a,b),$$ where the unions are over all
distinct $a$ in $\alpha $ and $|\alpha |\geq 2$ resp.\ distinct
pairs $a,b$ of elements in $\alpha $ and $|\alpha |\geq 3$.

\vspace{0.3cm}

\noindent \textbf{Theorem 2.3} \ \textit{Let $\sigma $ be a
geometric set of standard links. Then the sequence of
$\mathcal{R}$-modules: $$\CD \bigoplus \limits_{\alpha \in
\hat{\mathfrak{b}}(M)}\mathcal{R}
\mathcal{S}_{\mathcal{R}}(\alpha) @>{\chi _{\mathcal{R}}=\sum
\chi _{\mathcal{R},\alpha }}>>\mathcal{R}\hat{\mathfrak{b}}(M)
@>{\sigma }>> \mathcal{H}(M)\rightarrow 0.
\endCD $$
is exact for all choices of $\chi _{\mathcal{R}}$ determined by the
$\Delta$/$\Theta$-construction.}

\vspace{0.2cm}

\noindent \textsl{Proof.} \ Let $C=\mathfrak{b}(M)$ in theorem
2.1. Let $u$ act on the set $\mathfrak{b}(M)$ by defining $u\langle
\alpha_1,\ldots ,\alpha_r\rangle:=\langle 1,\alpha_1,\ldots
,\alpha_r\rangle$. Let $u^i$ act by composition for $i>1$ and let
$u^0:=id$. It follows that $\{u^{i}\langle \ \rangle |i\geq
0\}=\mathfrak{t}(M)$, and $u(\mathfrak{b}_r(M))\subset
\mathfrak{b}_{r+1}(M)$. The action defines the $R$-module
isomorphism $\mathcal{R}\hat{\mathfrak{b}}(M)\cong
R\mathfrak{b}(M)$. The set of relations in theorem 2.1 is indexed
by the union of the sets $$S(\alpha )=\bigcup_{a}T(a)\cup
\bigcup_{(a,b)}D(a,b),$$ where the unions are over all distinct
$a$ in $\alpha $ with $\alpha \setminus  a^{\circ} \notin
\mathfrak{t}(M)$ resp.\ all distinct pairs of
\textit{non-trivial} elements $(a,b)$ in $\alpha $ with $\alpha
\setminus \langle a^{\circ}, b^{\circ} \rangle \notin
\mathfrak{t}(M)$.

Consider the contributions from $T(1)$. These can be non-trivial only if
$\alpha $ also contains a
non-trivial conjugacy class. By definition $T(1)$ is a
set of generators of the group $\pi_1(LM,f_1)$ for a null
homotopic map $f_1: S^1\rightarrow M$, which represents a trivial
component of $K_{\alpha }$. Consider a corresponding self
homotopy of this trivial component of $K_{\alpha }$ to itself and
the resulting map $f: S^1\times S^1\rightarrow M$ with $f|S^1
\times *=f_1$. We can change the map by homotopy (compare section
6) such that $S^1\vee S^1\rightarrow M$ maps into a circle (the
image of $f(*\times S^1)$) without intersecting other components
of $K_{\alpha }$. There exists an extension to a map of torus into
that circle. The corresponding self homotopy of $K_{\alpha }$ will
give no contribution since the singular torus can be assumed not
intersecting with any other components of $K_{\alpha}$ by
transversality. Now $\pi_2(M)$ acts on the set $\pi_1(LM,f_1)$ by
changing maps on the top cell. So we can assume that
$T(1)=\mathfrak{S}$ without changing the image of $\chi $. Consider the decomposition $$S(\alpha
)=U(\alpha)\cup \bigcup_{a\neq
1}T(a)\cup
\bigcup_{(a,b)}D(a,b),$$ where
$U(\alpha )$ is equal to $T(1)$ if
$\alpha $ contains both a trivial conjugacy class and a
non-trivial conjugacy class, and is empty otherwise.

Now reorder the set $\cup_{\alpha}\mathcal{S}(\alpha )$ into $\cup_{\alpha}\mathcal{S}_{\mathcal{R}}'(\alpha )$, where
$$\mathcal{S}_{\mathcal{R}}'(\alpha ):=U(u\alpha )\cup
\bigcup_{a\neq 1}T(a)\cup \bigcup_{(a,b)}D(a,b).$$

\noindent The unions are restricted by the same conditions on $\alpha $ and $a,b$ resp.\ $a$ as before. Let the homomorphism $\chi _{\mathcal{R}}
=\sum_{\alpha}\chi _{\mathcal{R},\alpha}$ be defined on $\cup_{\alpha}\mathcal{S}_{\mathcal{R}}'(\alpha )$
by the $\Delta$/$\Theta$-construction as before with values now in $S\mathcal{R}\hat{\pi }(M)$.

The module $\oplus_{\alpha \in
\mathfrak{b}(M)}R\mathcal{S}_{\mathcal{R}}'(\alpha )$ admits the canonical
action of $u$ by shift of basis elements. Note that the
non-triviality restrictions in the definition of
$\mathcal{S}_{\mathcal{R}}'(\alpha)$ are preserved under multiplication by
$u$. Also for $\alpha \in \hat{\mathfrak{b}}(M)$ we have that
$\mathcal{S}_{\mathcal{R}}(\alpha )=\mathcal{S}_{\mathcal{R}}'(\alpha )$.

The main point is that the homomorphism $\chi _{\mathcal{R}}$ is compatible with ths actions,
i.e.
$\chi _{\mathcal{R},u\alpha }=u\chi _{\mathcal{R}, \alpha }$
for all $\alpha $. This follows because of the definition of $\chi _{\mathcal{R}}$ and condition 2 in the definition of geometric standard links.
So we can replace the free $R$-module with basis
$\cup_{\alpha \in \mathfrak{b}(M)}\mathcal{S}_{\mathcal{R}}'(\alpha )$
by the free $\mathcal{R}$-module with basis
$\cup_{\alpha \in \hat{\mathfrak{b}}(M)}
\mathcal{S}_{\mathcal{R}}(\alpha )$.  Obviously the restrictions that $\alpha \setminus  a^{\circ} $ resp.\
$\alpha
\setminus \langle a_1^{\circ}, a_2^{\circ} \rangle $ are not in
$\mathfrak{t}(M)$ for $\alpha \in \mathfrak{b}(M)$
translate into
$|\alpha |\geq 2$ resp.\ $|\alpha |\geq
3$ for $\alpha \in \hat{\mathfrak{b}}(M)$.
$\square$

\vspace{0.2cm}

\noindent \textbf{Remark.} \ The proof shows that there are
analogous presentations for $\mathcal{R}$-modules
$\mathcal{H}_{C^{\ast}}(M)$ and skein closed sets $C$.

\vspace{0.2cm}

For later purpose we state the following simple result at this
point.

\vspace{0.2cm}

\noindent \textbf{Lemma 2.4.} \textit{The natural map
$\mathfrak{H}(M)\rightarrow \mathcal{H}(M)$ induces the injective
homomorphism of abelian groups
$$\mathbb{Z}\mathfrak{H}(M)\rightarrow \mathcal{H}(M).$$}

\vspace{-0.3cm}

\noindent \textsl{Proof.} \ Consider the composition with the homomorphism
$\mathcal{H}(M)\rightarrow S\mathbb{Z}\hat{\pi }(M)$ induced by the
epimorphism $R\rightarrow \mathbb{Z}$, which is defined by
$q\mapsto 1, z\mapsto 0$ (compare section 1). This composition is injective. $\square$

\vspace{0.5cm}

\centerline{\textsc{3.\ HOMOTOPIES OF PATHS IN MAPPING SPACES}}

\vspace{0.3cm}

First we describe simple standard ideas from homotopy theory, which in fact are the basic ingredients for 2.1.

\vspace{0.2cm}

For $r\geq 1$ let $\tilde{\mathcal{M}}_r$ denote the space of
differentiable maps $P=\coprod_{r}S^1\rightarrow M$. The space
$\tilde{\mathcal{M}}_r$ is equipped with a topology in the
following way: A map $X\rightarrow \tilde{\mathcal{M}}_r$ is
continuous if the corresponding map $X\times P\rightarrow M$ is
continuous. Let $\mathcal{M}_r$ denote the quotient of this
mapping space by the permutation action of the symmetric group.
Let $\mathcal{M}=\cup_{r\geq 1}\mathcal{M}_r$. A set of
parametrized unordered standard links $K_{\alpha }$ can be
considered as a set of base points in the components of this
space. The \textit{link homotopy space} $\mathcal{M}_{\ell h}$ is
the subspace of those maps in $\mathcal{M}$ without intersections
of different components.

\vspace{0.2cm}

For each $\alpha $ there is the exact homotopy sequence of
\textit{pointed} sets: $$\CD \pi_1(\mathcal{M}_{\ell
h},K_{\alpha})\rightarrow \pi_1(\mathcal{M},K_{\alpha})\rightarrow
\pi_1(\mathcal{M},\mathcal{M}_{\ell h},K_{\alpha })@>\partial>>
\mathfrak{H}_{\alpha}(M)\endCD
$$
with $\partial $ surjective. By definition,
$\pi_1(\mathcal{M},\mathcal{M}_{\ell h},K_{\alpha }) $ is the set
of homotopy classes of paths in $\mathcal{M}$, which start at a
point of $\mathcal{M}_{\ell h}$ (i.e. a link map) and end in the
standard link $K_{\alpha }$. The action of
$\pi_1(\mathcal{M},K_{\alpha})$ on
$\pi_1(\mathcal{M},\mathcal{M}_{\ell h},K_{\alpha })$ is
transitive on the fibres $\partial^{-1}([f])$ for all $[f]\in
\mathfrak{H}_{\alpha}(M)$($[f]$ is the link homotopy class of a
link map $f$). If $\gamma_0$ is a loop in $K_{\alpha }$ and
$\gamma $ is a path ending in $K_{\alpha }$ then the action is
defined by juxtaposition. This defines a 1-1 correspondence
between the orbit space and $\mathfrak{H}_{\alpha }(M)$ and so
the $R$-isomorphism (for each $C\subset \mathfrak{b}(M)$)
$$\CD
\bigoplus \limits_{\alpha \in
\bar{C}}R(\pi_1(\mathcal{M},\mathcal{M}_{\ell h},K_{\alpha
})/\pi_1(\mathcal{M},K_{\alpha}))@>{\bar{\partial}}>>R\mathfrak{H}_{\bar{C}}(M).\endCD$$
\noindent By pull-back of the skein submodule $U_C(M)$ of skein
triples of links with wrapping invariant in $C$ we derive the
isomorphism:
$$\CD
\bigoplus \limits_{\alpha \in
\bar{C}}R(\pi_1(\mathcal{M},\mathcal{M}_{\ell h},K_{\alpha
})/\pi_1(\mathcal{M},K_{\alpha}))/(\bar{\partial}^{-1}(U_C))
@>{\partial'}>> \mathcal{H}_{C}(M).\endCD$$ The relation sets
$\bar{\partial}^{-1}(U_C)$ are generated by all
$q^{-1}\gamma_+-q\gamma_--z\gamma_0$. Here $\gamma $ is a
\textit{homotopy class} of paths in $\mathcal{M}$ (representing
an element in $\pi_1(\mathcal{M},\mathcal{M}_{\ell h},K_{\alpha
}$) up to the action of $\pi_1(\mathcal{M},K_{\alpha})$. From the
definition, $\gamma_{\pm}$ resp.\ $\gamma_0$ are \textit{any}
paths from $K_{\pm}$ resp.\ $K_0$ to $K_{\alpha }$ resp.\
$K_{\alpha '}$ for any skein triples $K_{\pm},K_0$ in $M$
($\omega (K_0)=\alpha '$). In section 5 we will use $\partial'$
to deduce the presentation 2.1 for the skein modules
$\mathfrak{H}_C(M)$. For this we need further notions and results
about homotopies of paths in the mapping spaces. The following is
due to X.\ S.\ Lin [L].

\vspace{0.2cm}

\noindent \textbf{Definition.} Let $j: P=\coprod_rS^1\rightarrow M$ be a differentiable
map. \newline
\textbf{a)} $j$ is called $1$-generic if it is an immersion with at most a single doublepoint,
which is \textit{rigid vertex}(this means that
the tangent vectors in the double point span a two-dimensional oriented subspace of the
tangent space of $M$ in that point).  \newline
\textbf{b)} \ $j$ is called \textit{$2$-generic} if it is
either an immersion with at most two rigid vertex doublepoints or it is an embedding except in
a single point, where the derivative vanishes.

\vspace{0.2cm}

Let $\tilde{\mathcal{M}}_1$ resp.\ $\tilde{\mathcal{M}}_2$ denote
the subspaces of $\tilde{\mathcal{M}}$ of $1$-generic resp.\
$2$-generic differentiable maps. Note that
$\tilde{\mathcal{M}}_1\subset \tilde{\mathcal{M}}_2$. Let
$\tilde{\mathcal{M}}_0$ denote the subspace of differentiable
embeddings. We consider maps $j_i: V^i\rightarrow
\tilde{\mathcal{M}}_i$ for compact oriented $i$-dimensional
manifolds $V_i$ and $i=1,2$. The singularity set $\mathfrak{S}_i$
of such a map is defined by
$\mathfrak{S}_i:=j_i^{-1}(\tilde{\mathcal{M}}_i\setminus
\tilde{\mathcal{M}}_0)$.

\vspace{0.2cm}

\noindent \textbf{Definition.} Let $i\in \{1,2\}$. A map $j_i
:V^i\rightarrow \tilde{\mathcal{M}}_i$ is called in
\textit{almost general position} if \newline \textbf{(i)}
$\mathfrak{S}_1$ is a finite collection of points in
the interior of $V^1$, resp.\ \newline \textbf{(ii)}
$\mathfrak{S}_2$ is a compact immersed $1$-manifolds with
transverse double points and $\mathfrak{S}_2\cap \partial V^2\subset
 j_2^{-1}(\tilde{\mathcal{M}}_1\setminus
\tilde{\mathcal{M}}_0)$, $j_2$  maps double points of $\mathfrak{S}_2$
to immersions with two rigid vertex double points. Each
boundary point of $\mathfrak{S}_2$ in the interior of $V^2$ is
mapped to a non-immersion in $\tilde{\mathcal{M}}_2$. All
embedded points, or boundary points of $\mathfrak{S}_2$ in
$\partial V^2$, are mapped into $\tilde{\mathcal{M}}_1$.

\vspace{0.2cm}

Let $\mathfrak{S}'\subset \mathfrak{S}_2$ denote the set of non-embedding points
or boundary points in the interior of $V^2$.

\vspace{0.2cm}

\noindent \textbf{Remarks.} \newline \textbf{a)} \ The restriction
of a map $V^2\rightarrow \tilde{\mathcal{M}}_2$ in almost general
position to the boundary is a map $\partial V^2\rightarrow
\tilde{\mathcal{M}}_1$ in almost general position. \newline
\textbf{b)} \ All differentiable maps in the image of a component
of $V^i\setminus \mathfrak{S}_i$ are isotopic embeddings,
$i=1,2$. \newline \textbf{c)} \ All differentiable maps in the
image of a component of  $\mathfrak{S}_2\setminus \mathfrak{S}'$
are differentiable isotopic immersions with a single rigid vertex
double point.\newline  \textbf{d)} \ A neighbourhood of an
interior boundary point of $\mathfrak{S}_2$ in $\mathfrak{S}_2$
describes the shrinking of a kink in a component to a point  (see [L]).

\vspace{0.2cm}

\noindent \textbf{Theorem 3.1. (Lin [L])} \ \textit{Let $j_i:
V^i\rightarrow \tilde{\mathcal{M}}$ be a map, whose restriction
to the bounday is in almost general position (in particular maps
into $\tilde{\mathcal{M}}_{i-1}$). Then $j_i$ can be approximated
rel boundary by a map $j_i' :V^i\rightarrow
\tilde{\mathcal{M}}_i$ in almost general position.}

\vspace{0.2cm}

We are only interested in that part of the singularity sets
$\mathfrak{S}_i$, which maps to immersions with only double
points of different components. This defines submanifolds
$\mathfrak{D}_i\subset \mathfrak{S}_i\subset V_i$. It follows
from Lin`s work that the immersed singularity manifold
$\mathfrak{D}_2$ satisfies the following properties.

\vspace{0.2cm}

\noindent \textbf{Proposition 3.2.} \textit{Let $j_2:
V^2\rightarrow \tilde{\mathcal{M}}_2$ be a map in almost general
position. Then the following holds:}
\newline
\textbf{a)} \ \textit{The set $\mathfrak{D}_2$ is a properly
immersed $1$-dimensional manifold in $V^2$.}\newline \textbf{b)}
\ \textit{$j_2$ maps all points in a component of $V^2\setminus
\mathfrak{D}_2$ to link homotopic link maps (i.e.\ without
singularities of different components, but with possible
singularities of type $\tilde{\mathcal{M}}_2$ of same
components)}.\newline \textbf{c)} \ \textit{If $\mathfrak{D}'$
denotes the double point set of $\mathfrak{D}_2$ then all points
in a component of $\mathfrak{D}_2\setminus \mathfrak{D}'$ are
mapped to immersions with a single double point of different
components and possibly a further double point of same
components.  $\square$}

\vspace{0.2cm}

\noindent \textbf{Definition.} \ Let $j_2: V^2\rightarrow
\tilde{\mathcal{M}}_2$ be in almost general position. A point in
$\mathfrak{D}'$ maps to an immersion with two double points of
different components. We call such a point \textit{paired} if
both double points are intersection points of the \text{same}
pair of components. Otherwise the point in $\mathfrak{D}'$ is
called \textit{unpaired}.

\vspace{0.2cm}

A map $V^i\rightarrow \mathcal{M}$ is \textit{in almost general
position} if it lifts to a map $V^i\rightarrow
\tilde{\mathcal{M}}$ in almost general position. The definition
of singularity sets $\mathfrak{D}_i$ and $\mathfrak{D}'$
obviously applies to all maps $V^i\rightarrow \mathcal{M}$, which
can be lifted.

\vspace{0.2cm}

Let $\mathcal{M}^{\ast}$ denote the set of those maps
$P\rightarrow M$ with a single rigid vertex intersection point of
different components and possible arbitrary self-intersections of
components disjoint from the intersection point. Let
$\mathfrak{H}^1(M):=\pi_0(\mathcal{M}^{\ast})$. We can obviously
replace $\mathcal{M}^{\ast}$ by $\mathcal{M}^{\ast}\cap
\mathcal{M}_2$ without changing the set of path components.

\vspace{0.2cm}

We prove a few results about link homotopy and skein theory
related to the ideas above.

\vspace{0.2cm}

\noindent \textbf{Proposition.} \ \textit{Let $j: M\hookrightarrow N$ be a inclusion of
oriented $3$-manifolds. Assume that $j$ induces a bijection of link homotopy
sets
$j^0: \mathfrak{H}(M)\rightarrow \mathfrak{H}(N)$
and a surjection
$j^1: \mathfrak{H}^1(M)\twoheadrightarrow \mathfrak{H}^1(N)$.
Then the induced homomorphism
$j_*: \mathcal{H}(M)\rightarrow \mathcal{H}(N)$
is an isomorphism.}

\vspace{0.2cm}

\noindent \textsl{Proof.} \ The isomorphism $j_* : R\mathfrak{H}(M)\rightarrow
R\mathfrak{H}(N)$ descends to an epimorphism $\mathcal{H}(M)\rightarrow \mathcal{H}(N)$,
which is injective if $U_N=j_*(U_M)$ (note that $j_*(U_M)\subset U_N$ always holds). But
this follows from the surjectivity of $j^1$. $\square$

\vspace{0.2cm}

\noindent \textbf{Corollary.} \ \textit{Assume that $j_*: \mathfrak{H}(M)\rightarrow
\mathfrak{H}(N)$ is bijective and $j$ induces an isomorphism of fundamental groups.
Then the induced homomorphism
of $q$-homotopy skein modules is an isomorphism. $\square$}

\vspace{0.2cm}

\noindent \textbf{Examples.} \ The inclusion $M\hookrightarrow
\hat{M}$ induces the isomorphism $\mathcal{H}(M)\rightarrow
\mathcal{H}(N)$. Let $e\subset M$ be a (possibly) fake open
$3$-cell. It is known that $M\setminus e\hookrightarrow M$
induces the bijection $\mathfrak{H}(M\setminus e)\rightarrow
\mathfrak{H}(M)$ for $M=S^3$ (see [H]). The corresponding result
for $M\neq S^3$ is not known (the surjectivity part is obvious).
In section 8 we will show that $\mathcal{H}(M\setminus
e)\rightarrow \mathcal{H}(M)$ is an isomorphism.

\newpage

\centerline{\textsc{4.\ THE BASIC CONSTRUCTION}}

\vspace{0.3cm}

For $\alpha \in \mathfrak{b}(M)$ let $\mathcal{P}(\alpha )$
denote the set of paths in almost general position $\gamma :
I\rightarrow \mathcal{M}_1\subset \mathcal{M}$ in the component
of $\mathcal{M}$ determined by $\alpha $. In particular, $\gamma
\in \mathcal{P}(\alpha )$ implies $\gamma(0),\gamma(1)$ are
embeddings and $\gamma (t)$ has free homotopy classes given by
$\alpha $ and is not an embedding for only a finite number of
points.

\vspace{0.2cm}

\noindent \textbf{Definitions 4.1.} \ Let $t_1,\ldots ,t_k$ be the
parameters for which $\gamma (t_i)$ is an immersion with a single
intersection point of \textit{different} components (the
\textit{singularity parameters for link homotopy}). Let $K_{t_i}$
be the result of smoothing the corresponding rigid vertex
immersion in the usual way. We define $\varepsilon_i\in \{\pm
1\}$ by $\varepsilon_i=\pm 1$ if $K_{t_i-\delta}=K_{\pm}$ and
$K_{t_i+\delta}=K_{\mp}$ for $\delta
>0$ small enough ($K_{\pm}$ are the links defined by the
double point). Let $\varepsilon (\gamma ):=
\sum_{i=1}^k\varepsilon_i$ be the \textit{index} and let
$\sum_i^k|\varepsilon_i|$ be the \textit{length} of $\gamma $.

\vspace{0.2cm}

\noindent \textbf{Remark.} \ By 3.2 the index of paths in almost
general position is invariant under homotopies in $\mathcal{M}$
with endpoints in $\mathcal{M}_{\ell h}$. This can be shown by
lifting to a homotopy in the covering space $\tilde{\mathcal{M}}$
and applying 3.1.

\vspace{0.2cm}

For $\alpha \in \mathfrak{b}(M)$ let $\lfloor \alpha \rfloor \subset
\mathfrak{b}(M)$ denote the set of all first order descendants of $\alpha $
(compare section 2).

\vspace{0.2cm}

\noindent \textbf{Definition.} \ For each \textit{elementary}
path $\gamma $ with a single singularity point $t_1$ of index
$\varepsilon_1 $ define
$$s: \mathcal{P}(\alpha )
\longrightarrow \mathbb{Z}[q^{\pm 1}] \mathfrak{H}_{\lfloor \alpha
\rfloor}(M).$$ by $s(\gamma ):=
\varepsilon_1q^{\varepsilon_1}K_{t_1}$. Extend the definition by
the juxtaposition formula (assuming $\gamma_1(1)=\gamma_2(0)$)
$$s(\gamma_1 \gamma_2):=s(\gamma_1)+q^{2\varepsilon
(\gamma_1)}s(\gamma _2),$$ where $\gamma_1\gamma_2$ is the usual
product of paths.

\vspace{0.2cm}

We need some elementary properties of the map $s$. Throughout, if not indicated otherwise,
\textit{path} means path in almost general position. Note that  $s$ does not
change if we pre- or post compose by a path (in almost general
position) in $\mathcal{M}_{\ell h}$.

\vspace{0.2cm}

\noindent  \textbf{Lemma 4.2.} \ \textit{For $\gamma $ a path with
singularities $0<t_1< \ldots <t_k<1$ of indices
$\varepsilon_1,\ldots ,\varepsilon_r$ the following holds: $$\CD
s(\gamma )= \varepsilon_1q^{\varepsilon_1}K_{t_1}
+\varepsilon_2q^{2\varepsilon_1+\varepsilon_2}K_{t_2}
+\varepsilon_3q^{2(\varepsilon_1+\varepsilon_2)+\varepsilon_3}K_{t_3}
\\ +\ldots + \varepsilon_kq^{2(\varepsilon_1+\ldots
+\varepsilon_{k-1})+\varepsilon_k}K_{t_k}.\endCD$$ }

\noindent \textbf{Remark.} If $\gamma $ has non-trivial
index then $s(\gamma )\neq 0$. Also the image of
$s(\gamma )$ in $\mathcal{H}(M)$ cannot be trivial. In order to
show this use 2.4. to see that $\varepsilon (\gamma )$ can be computed from
the image in the homotopy skein module.

\vspace{0.2cm}

\noindent \textbf{Lemma 4.3.} \ \textit{Let $\gamma ^{-1}$ the inverse path of $\gamma $.
Then
$$s(\gamma \gamma ^{-1})=s(\gamma ^{-1}\gamma )=0.$$}
\noindent \textsl{Proof.} \ This is proved by induction on the length of paths. $\square $

\vspace{0.2cm}

The following is checked by direct computation.

\vspace{0.2cm}

\noindent \textbf{Lemma 4.4.}\newline \textbf{a)} \  \textit{For each path $\gamma $:
$s(\gamma^{-1})=-q^{-2\varepsilon (\gamma )}s(\gamma )$}\newline \textbf{b)} \
\textit{For two paths $\gamma _1,\gamma_2$ with the same initial and end points we
have $$s(\gamma _2)-s(\gamma _1)=q^{2\varepsilon (\gamma_1)}s(\gamma_1^{-1}\gamma
_2)$$ \textbf{c)} Let $\gamma_1,\gamma_2$ be paths with $\gamma_1(1)=\gamma_2(0)$,
then $$s(\gamma_1\gamma_2)-s(\gamma_2\gamma_1)=(1-q^{\varepsilon
(\gamma_2)})s(\gamma_1)-(1-q^{\varepsilon (\gamma_1)})s(\gamma_2) \quad  \square  $$ }

\vspace{-0.3cm}

The map $s$ is not always
\textit{commutative} with respect to composition. It is
commutative, if $\varepsilon (\gamma_1)=\varepsilon (\gamma_2)=0$.

\vspace{0.2cm}

An important consequence of 4.4 c) is the following result:

\vspace{0.2cm}

\noindent \textbf{Corollary 4.5.} \ \textit{ Let $\gamma_1$ be a
loop in $\gamma_0(1)$. Then the following holds:
$$s(\gamma_0\gamma_1\gamma_0^{-1})=(1-q^{2\varepsilon
(\gamma_1)})s(\gamma_0)+q^{2\varepsilon (\gamma_0)}s(\gamma_1). \quad \square $$}

\vspace{-0.3cm}

If $\varepsilon (\gamma_1)=0$ then conjugation by some
path $\gamma_0$ (this is changing the base point) just amounts to multiplication
by $q^{2\varepsilon (\gamma_0)}$. In particular this applies when the base point
on a loop with trivial index is moved on the loop.

\vspace{0.2cm}

The map $s$ measures the \textit{skein difference} between
initial and end point of a path. In order to define absolute
skein invariants we need to include the links  $\gamma (0)$ and
$\gamma (1)$  into the definition.

\vspace{0.2cm}

\noindent \textbf{Definitions 4.6.} \ Let $\gamma \in \mathcal{P}(\alpha )$. Then we define
$$s_f(\gamma)=q^{2\varepsilon (\gamma)}\gamma (1)+zs(\gamma )
\in R\mathfrak{H}_{\bar{\alpha }}(M)$$
and
$$s_{if}(\gamma )=q^{2\varepsilon (\gamma )}\gamma (1)-\gamma (0)
+zs(\gamma )\in R\mathfrak{H}_{\bar{\alpha }}(M).$$

\vspace{0.2cm}

Let $\mathcal{P}_{\alpha }$ be the set of paths in almost general
position with end point in $K_{\alpha }$. Recall that
$\pi_1(\mathcal{M},\mathcal{M}_{\ell h},K_{\alpha })$ is the set
of homotopy classes of paths with end point $K_{\alpha }$ and
initial point in $\mathcal{M}_{\ell h}$.

\vspace{0.2cm}

The map
$$s_f :\mathcal{P}_{\alpha }\longrightarrow R(\mathfrak{H}_{{\lfloor \alpha
\rfloor}}(M) \cup \{K_{\alpha}\})$$ assigns to each link with
wrapping invariant $\alpha $ a linear combination of $K_{\alpha }$
and links with fewer components. It describes the skein expansion
of a link along a path $\gamma $ in terms of a standard link and
links with fewer components.

\vspace{0.2cm}

The properties of the map $s$ translate to important properties of
$s_f$ and $s_{if}$.

\vspace{0.2cm}

\noindent \textbf{Proposition 4.7.} Let $\gamma_1,\gamma_2$ with $\gamma_1(1)=
\gamma_2(0)$. Then
$$s_{if}(\gamma_1\gamma_2)=s_{if}(\gamma_1)+q^{2\varepsilon
(\gamma_1)}s_{if}(\gamma_2), \ \text{and} $$
$$s_f(\gamma_1\gamma_2)=s_f(\gamma_1)+q^{2\varepsilon (\gamma_1)}s_{if}(\gamma_2).$$
\noindent \textsl{Proof.} \ Both results are immediate by
computation. For the second formula use that $s_f(\gamma
)=s_{if}(\gamma )+\gamma (0)$. $\square $.

\vspace{0.2cm}

\noindent \textbf{Remark 4.8.} \newline
\textbf{a)} \ The second formula generalizes the skein relations
$K_+=q^2K_-+zqK_0$ and $K_-=q^{-2}K_+-zq^{-1}K_0$ in the following way:
Let $\gamma_1$ be an elementary path with singularity parameter $t_1$ of index $\varepsilon
(\gamma_1)=\varepsilon_1$. Then it follows from the above
\begin{eqnarray*}
s_f(\gamma_1\gamma_2)=&s_f(\gamma_1)+q^{2\varepsilon
(\gamma_1)}(s_f(\gamma_2)-\gamma_2(0))\\
\ =&q^{2\varepsilon (\gamma_1)}s_f(\gamma_2)+q^{2\varepsilon (\gamma_1)}
(\gamma_1(1)-\gamma_2(0))+zs(\gamma_1) \\
\ =&q^{2\varepsilon_1}s_f(\gamma_2)+\varepsilon_1zq^{\varepsilon_1}K_{t_1}
\end{eqnarray*}

The inital point of $\gamma_1\gamma_2$ is a link $K_{\pm}$, the
initial point of $\gamma_2$ is $K_{\mp}$ and the smoothing $K_0$
is $K_{t_1}$. \newline \textbf{b)} \ The result holds in
particular for all compositions with
$\gamma_1(1)=\gamma_2(0)=\gamma_2(1)$. Thus it describes the
action of the loops in $K_{\alpha }$ on paths joining links with
$K_{\alpha }$.

\vspace{0.2cm}

To combine the results of this section with section 3 define the homomorphism
$$\CD R\bar{C}@>{\sigma '}>>
\bigoplus \limits_{\alpha \in
\bar{C}}R(\pi_1(\mathcal{M},\mathcal{M}_{\ell h},K_{\alpha
})/\pi_1(\mathcal{M},K_{\alpha}))/(\bar{\partial
}^{-1}(U_C)),\endCD$$ by mapping the basis element determined by
$\alpha \in \bar{C}$ to the constant loop in $K_{\alpha }$.

\vspace{0.2cm}

\noindent \textbf{Proposition.} \textit{The homomorphism
$\sigma'$ above is surjective.}

\vspace{0.2cm}

\noindent \textsl{Proof.} \  Represent elements of
$\pi_1(\mathcal{M},\mathcal{M}_{\ell h},L_{\alpha
})/\pi_1(\mathcal{M},K_{\alpha})$ by paths in almost general
position, which join a link $K$ to the basepoint $K_{\alpha }$.
Then proceed by induction on the complexity $\ell (\gamma )=
(|\alpha |, \text{length of}\ \gamma)$ (see 4.1). If there are no
singular parameters then $\gamma $ can be homotoped into the
constant map. Otherwise consider the first $t_1>0$ with
$\gamma(t_1)\in \mathfrak{D}_1$. Let $\gamma =\gamma_1\gamma_2$,
where $\gamma_1$ is an elementary path starting in $K$ and ending
in $K_{t_1+\delta }$ for $\delta
>0$ sufficiently small. Choose a path $\gamma_0$, which joins $K_{t_1}$ to its standard link.
The triple $\gamma ,\gamma_2,\gamma_0$ is a skein triple of paths. So we can replace
$\gamma $ by paths of smaller complexities $\gamma_2, \gamma_0$.$\square$

\vspace{0.5cm}

\centerline{\textsc{5.\ PROOF OF THEOREM 2.1.}}

\vspace{0.3cm}

We can assume that $C$ is a skein closed subset of
$\mathfrak{b}(M)$. The set $C$ decomposes
into the union $C=\cup_{i\geq 1}C(i)$, where $C(i)$ is the set of
all elements of $C$ of length $\leq i$.

\vspace{0.2cm}

The idea of the proof is to define a homomorphism
$$\rho : \mathcal{H}_C(M)\longrightarrow RC/W$$
for the submodule $W:=im(\chi )$. This will be done
such that
\begin{enumerate}
\item the composition
$$\CD
RC@>{\sigma}>>\mathcal{H}_C(M)@>{\rho }>>RC/W \endCD$$
is the natural projection (this proves $ker(\sigma )\subset W$),
and
\item $\sigma (W)=0$ (or $W\subset ker(\sigma )$ for $W\subset RC$).
\end{enumerate}
It follows that $ker(\sigma )=W$, which is the assertion of
theorem 2.1.

\vspace{0.2cm}

The homomorphism $\rho $ will be defined using the isomorphism $\partial '$ from section 3:
First define a skein invariant map on
$\mathcal{P}_{\alpha }$. Then prove that it only depends on the
link homotopy class of initial points of paths in $\mathcal{P}_{\alpha }$.

\vspace{0.2cm}

We will construct $\rho $ and $W$ inductively through the sequence of homomorphisms:
$$\rho_r: \mathcal{H}_{C(r)}(M)\longrightarrow A(r)$$
for the modules
$A(r):=RC(r)/W(r)$, where $W(r)$ is the submodule of $RC(r)$, which is generated by images
of $\chi_{\alpha }$ for
all $\alpha \in C(r)$. The result follows from the fact that $\mathcal{H}_C(M)$ is the direct
limit of the modules
$\mathcal{H}_{C(r)}(M)$ for $r\rightarrow \infty$, which is obvious from the definitions.

\vspace{0.2cm}

Define
$\rho_1(K)=\rho_1(\gamma ):=\alpha \in RC(1)$
for each knot $K$ with $\omega (K)=\alpha $ resp.\ path joining $K$ with
$K_{\alpha }$.
This is well defined and extends to the homomorphism
$$\mathcal{H}_{C(1)}(M)\longrightarrow RC(1)=:A(1).$$
There are no skein relations on this level, which need to be considered.

\vspace{0.2cm}

We assume inductively that we have defined
$$\rho_{r-1}: \mathcal{H}_{C(r-1)}(M)\longrightarrow A(r-1).$$
For links with wrapping invariant in $C(r-1)$ we define $\rho_r$ by composition of
$\rho_{r-1}$ with the natural
homomorphism: $$j_{r-1,r}: A(r-1)=RC(r-1)/W(r-1)\rightarrow RC(r)/W(r)=A(r).$$
We have to define $\rho_r$ on $\mathfrak{H}_{C(r)}(M)$.
Choose a path $\gamma $ in almost general position, which joins $K$ with $K_{\alpha}$ for $\alpha \in C(r)$. Consider
$$s_f(\gamma )=q^{2\varepsilon (\gamma )}K_{\alpha }+zs(\gamma )
\in R(\mathfrak{H}_{\lfloor \alpha \rfloor}(M) \cup \{K_{\alpha}\})
\subset R(\mathfrak{H}_{C(r-1)}(M)\cup \{K_{\alpha }\}).$$
Define the $R$-homomorphism $\mathfrak{p}(\alpha )$ :
$$R(\mathfrak{H}_{C(r-1)}(M)\cup \{K_{\alpha }\})
\rightarrow R(C(r-1)\cup \{\alpha \})/W(r-1)=:A'(r-1)$$
by composition with $\rho_{r-1}$ on $\mathfrak{H}_{C(r-1)}(M)$ (actually we first have to
project this into $\mathcal{H}_{C(r-1)}(M)$), and by mapping $K_{\alpha }$ to $\alpha $. This is
well defined since $W(r-1)$ does not contain any terms involving $K_{\alpha }$.

Finally let
$\rho '(\alpha ):=\mathfrak{p}(\alpha )\circ s_f :\mathcal{P}_{\alpha }\rightarrow A'(r-1)$.
The next step is to change the module $A'(r-1)$ to the module
$A(r)$ in such a way that the map $\rho '$ does not depend on the
choice of path but only on the initial point.

\vspace{0.2cm}

We have to consider the following two ways of changing $\gamma $:

\begin{enumerate}
\item The path $\gamma $ is changed in its homotopy class
in $\pi_1(\mathcal{M},\mathcal{M}_{\ell h},L_{\alpha })$
\item The path is changed by the action of $\pi_1(\mathcal{M},K_{\alpha})$.
\end{enumerate}

\vspace{0.2cm}

\noindent \textbf{Claim.} \ \textit{The contributions of $\rho '(\alpha )$ defined by} (1) \textit{resp.\ }(2) \textit{are contained in
$pr(\chi_{\alpha }(\mathcal{S}(\alpha))$, where $pr$ is the projection
$R(C(r-1)\cup \{\alpha \}) \rightarrow A'(r-1)$.}

\vspace{0.2cm}

Let $K\in \mathfrak{H}_{\alpha}(M)$
and $\gamma \in \mathcal{P}_{\alpha }$ with
$\gamma (0)=K$. Then define
$$\rho_r(K)=pr' \circ \mathfrak{p}(\alpha )\circ s_f(\gamma )$$
for some path $\gamma \in \mathcal{P}_{\alpha }$ with $\gamma (0)=K$. Here we use the
canonical projection
$pr' :A'(r-1)\rightarrow A(r)$.

\vspace{0.2cm}

\noindent \textbf{Remarks.} \newline \textbf{a)} \ The
contributions to $W$ from 1) resp.\ 2) are
\textit{local} resp.\ \textit{global} obstructions in the
construction of finite type link homotopy invariants of Homfly type
[Ka-L].
\newline \textbf{b)} \ If the independence of
the choice of paths has been established then skein invariance is
provided just by the very definition of the map $s_f$ and its
skein property stated in 4.8).
\newline \textbf{c)} Obviously
$\sigma (W)=0$ because $W$ is defined from terms $s_{if}(\gamma
)$ for loops in $K_{\alpha }$.

\vspace{0.2cm}

\noindent \textbf{i.) Discussion of the local contributions.}

\vspace{0.2cm}

Now a homotopy of some representative path in  $\pi_1(\mathcal{M},\mathcal{M}_{\ell
h},L_{\alpha })$ can be split up in homotopies, which move the
initial point in $\mathcal{M}_{\ell h}$, and homotopies
relative to the initial and end point.
Obviously $s_f$ is invariant in the first case. So we
only need to consider the second case. Let $\gamma '$ be a path
homotopic to $\gamma $, both paths in almost general position.
Consider the loop $\gamma^{-1}\circ \gamma '$ in $K_{\alpha }$.
This loop is null-homotopic in $\mathcal{M}$ thus it bounds a
disk in $\mathcal{M}$, which lifts  to $\tilde{\mathcal{M}}$. So
we can apply 3.1 and approximate the map of this disk into
$\tilde{\mathcal{M}}$ by an map in almost general position and
project this to $\mathcal{M}$. The singularity set
$\mathfrak{D}_2$ of this map of a disk (the set of points, for
which we have immersions with double points of different
components) is a properly embedded immersion by 3.2.

Now consider the deformation of $\gamma $ into $\gamma '$ defined
by the disk map in almost general position. Consider the disk
bounding the union of two intervals identified in their end
points. Homotope the two bounding intervals into each other across
the disk. By 3.2.c) the invariant $s(\gamma )$ can only change
when we homotope the path across a point of $\mathfrak{D}'$. In
fact note that in the homotopy pairs of intersection points with
$\mathfrak{D}_2$ can appear or vanish if we move the path over a
\textit{critical point} of $\mathfrak{D}_2$ with respect to a
suitable height function. But the two contributions here are
easily seen to cancel. We have to consider the situation in the
picture for $\varepsilon (\gamma_1)=\varepsilon (\gamma_2)\in
\{-2,0,+2 \}$. Up to cyclic permutation the indices of the
singular points on the circle are given by $(+1,+1,-1,-1)$.

\vspace{1.3cm} \setlength{\unitlength}{2mm}
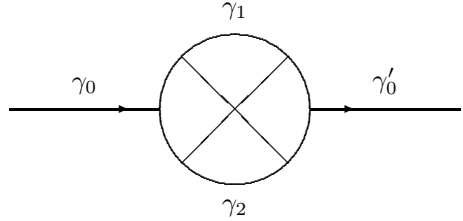
\begin{center}
\begin{picture}(30,1)
\drawline(0,0)(10,0) \put(0,0){\vector(1,0){8}}
\put(15,0){\arc(0,5){360}} \drawline(18.53,3.53)(11.46,-3.53)
\drawline(18.53,-3.53)(11.46,3.53) \drawline(20,0)(30,0)
\put(20,0){\vector(1,0){3}}
\put(5,1){\makebox(0,0)[b]{$\gamma_0$}}
\put(25,1){\makebox(0,0)[b]{$\gamma'_0$}}
\put(15,6){\makebox(0,0)[b]{$\gamma_1$}}
\put(15,-6){\makebox(0,0)[t]{$\gamma_2$}}
\end{picture}
\end{center}
\vspace{1.3cm}

Using 4.7. we compute the difference $s_f(\gamma_0
\gamma_1\gamma_0')-s_f(\gamma_0 \gamma_2\gamma_0')$:
$$s_f(\gamma_0\gamma_i\gamma_0')=s_f(\gamma_0)+q^{2\varepsilon
(\gamma_0)}s_{if}(\gamma_i)+q^{2\varepsilon (\gamma_0)+\varepsilon
(\gamma_i))}s_{if}(\gamma_0'), \  \text{so}$$

\vspace{-0.7cm}

\begin{eqnarray*}
\ & s_f(\gamma_0\gamma_2\gamma_0')-s_f(\gamma_0\gamma_1\gamma_0')
\ & = q^{2\varepsilon (\gamma_0)}(s_{if}(\gamma_2)-s_{if}(\gamma_1)) \\
= \ &  q^{2\varepsilon (\gamma_0)}(s_f(\gamma_2)-s_f(\gamma_1))
\ & = q^{2(\varepsilon (\gamma_0)+\varepsilon (\gamma_1))}s(\gamma_2^{-1}\gamma_1).
\end{eqnarray*}

It follows that we can focus on computing the value of
$s(\gamma_2^{-1}\gamma_1)$. By 4.5, a change of basepoint on the
loop alters the contribution only by multiplication with a power
of $q$.

Each point of $\mathfrak{D}'$ is an immersion with two double points of
different components. Thus we can compute the contribution $\Delta (L)$
from a loop in some link $K=K_{+-}$ (the two lower indices corresponding to the two double
points of the immersion $L$ in the crossing, the signs over the arrows indicating the indices):
$$\CD K_{+-}@>{+1}>>K_{--}@>{-1}>>K_{-+}@>{-1}>>K_{++}@>{+1}>>K_{+-}
\endCD$$
By 4.2. we compute:
$\Delta (L)=qK_{0-}-qK_{-0}-q^{-1}K_{0+}+q^{-1}K_{+0}$.

\vspace{0.2cm}

\noindent \textbf{Case 1.} \ Assume that the point in $\mathfrak{D}'$ is \textit{unpaired}
(compare 4.3). Then the $+$ resp.\ $-$ crossings in the formula are still crossings of different
components and the computation can be proceeded:
$$\Delta (L)=qK_{0-}-qK_{-0}-q^{-1}(q^2K_{0-}+qzK_{00})+q^{-1}(q^2K_{-0}+qzK_{00})=0.$$
\noindent \textbf{Case 2.} Assume that the point in $\mathfrak{D}'$ is \textit{paired}. Then we
can simplify:
$$\Delta (L)=z(q-q^{-1})(K_{0+}-K_{+0}).$$
Those are local contributions, which are not trivial in general. These are precisely the terms,
which are taken care of in the submodule $T'$. We have to show that $j_{r-1,r}\circ \rho_{r-1}(\Delta (L))=0$. Note that $\rho_{r-1}$ is well defined on the links $K_{0+},K_{+0}\in
\mathfrak{H}_{C(r-1)}(M)$.

Note that the links $K_{+0},K_{0+}$ result from $K=K_{++}$ by attaching bands to $K$, possibly
in different ways. Assume that the crossings, which are smoothed, are crossings of the first
two components of $K$.
Consider a path in almost general position joining $K$ and $K_{\alpha}$.
This path can be followed by a homotopy of the bands in the complement of the last $r-2$
components of $K$, with the end-points of the bands possibly moving in the first two
components. This follows from transversality. In fact all crossing changes of components can
be avoided just by moving the bands aside or the band is isotoped in $M$ with the link. So the
path above also provides two paths, which join $K_{+0}$ resp.\ $K_{0+}$ with two links
resulting from $K_{\alpha}$ by attaching bands. Consider the crossing changes for the two
deformations. This means applying e.g. relations
\begin{eqnarray*}
K_{0+,+}&=q^2K_{0+,-}+qzK_{0+,0}\\
K_{+0,+}&=q^2K_{+0,-}+qzK_{+0,0}
\end{eqnarray*}
(and similarly for $K_{0+,-}\mapsto
K_{0+,+}$). So the difference terms will still be a linear combination of terms
$K_{0+}'-K_{+0}'$ for links of $\leq r-2$ components resulting from the construction above on
links with $r-1$ components. Note that crossing changes of the first two components are not
contributing anything. So the additional contributions are already taken care of in the
induction.

Finally we have to show that the contribution $\Delta (L)$ does not depend on the explicit
choice of bands. Consider a homotopy of the bands relative to the attaching arcs. Of course
crossings of the band with any of the first two components do not contribute. Consider the
crossing of a band with any of the last $r-2$ components. The resulting difference in the
expansion of $K_{0+}-K_{+0}$ is easily computed. In fact, at such a crossing we have to
consider for e.g.\ $K':=K_{0+}$ the change from
$K_{+-}'$ to $K_{-+}'$ (these are different crossings now).

We compute:
$$K_{+-}'=q^2K_{--}'+qzK_{0-}'=K_{-+}'+qz(K_{0-}'-K_{-0}'),$$
where again we have smoothed a pair of components. The differnce just appears in the form
$$z(q-q^{-1})qz(K_{0+}'-K_{+0}')$$
in $\Delta (L)$. Since $K'=K_{0+}$ has $r-1$ components this term is already contained in the
submodule $W(r-1)$. The same argument works for
$K_{+0}$. Thus we see that the band can be homotoped freely in any position with endpoints on
the first two components. So only this homotopy class has to be considered. Also any kind of
twisting or knotting or linking of the band with other components can be neglected by the
arguments above.

We only have to account \textit{differences} of expansions in terms of standard links
resulting from attaching two bands to a standard link. So it suffices to consider all
differences with respect to a fixed choice determined by the basing. Recall that the basing of the
two first components determines a band (which runs from the first component to the basepoint
along the basing arc of the first component, then along the basing arc of the seond component
from the basepoint to the second component). The basing is changed by composition with loops in $*$.
In this way all possible bands (up to homotopy) appear. Moreover, we can slide
the band along the first component resp.\ along the second component. So we can change the
group element by multiplication with powers of $a$ from the left and powers of $b$ from
the right. For a more detailed consideration of basings and band operations see [K3]. So the $\Delta$-construction  generates the module $W'(r)$ of $W(r)$ (modulo
lower order constributions) of local contributions \newline
\noindent $z(q-q^{-1})(K_{0+}-K_{+0})$.

\vspace{0.2cm}

\noindent \textbf{ii) Discussion of global contributions.}

\vspace{0.2cm}

Up to this point we have defined a map:
$$\mathcal{P}_{\alpha }
\rightarrow RC(r)/(W(r-1)+W'(r)),$$ which now only depends on the
image of $\gamma \in \mathcal{P}_{\alpha} $ in
$\pi_1(\mathcal{M},\mathcal{M}_{\ell h},K_{\alpha })$. In order
to define the map $\rho_r$ we consider the change of $s_f(\gamma
)$ by the action of $\pi_1(\mathcal{M},K_{\alpha})$.

By 4.7, for $\gamma ,\gamma '\in \mathcal{P}_{\alpha }$ with
$\gamma '(0)=K_{\alpha }$:
$$s_f(\gamma \gamma ')=s_f(\gamma )+q^{2\varepsilon (\gamma )}s_{if}(\gamma ')$$
So let $W''(r)$ (mod $W(r-1)+W'(r)$) denote the submodule generated by the images
$s_{if}(\gamma ')$ for all loops in $K_{\alpha }$ in almost general position.

\vspace{0.2cm}

Then theorem 2.1 follows if we can prove that $W''(r)$ is generated by the $\Theta$-construction.

\vspace{0.2cm}

By 4.7 we only have to consider a generating set of
$\pi_1(\mathcal{M},K_{\alpha})$. Recall that
$\pi_1(\mathcal{M},K_{\alpha})$ is a subgroup of
$\pi_1(\tilde{\mathcal{M}},K_{\alpha}')$ for some lift
$K_{\alpha}'$ into the covering space. There are loops in
$\mathcal{M}$, which do not lift to a loop in the covering space.
But such a lift is a path joining two different orderings
$K_{\alpha}'$ and $K_{\alpha}''$ of $K_{\alpha }$. Obviously the
corresponding permutation of components of $K_{\alpha }$ can only
permute components with the same homotopy classes.
Now recall condition 3 in the definition of a geometric standard
basis. We can always compose with a loop of embeddings in
$\mathcal{M}_0\subset \mathcal{M}$, which permutes the components
in arbitrary possible ways, and gives no contributions to $s_f$
(use 4.7). Thus, in order to compute $W''(r)$, we can restrict to
the subgroup $\pi_1(\tilde{\mathcal{M}},K_{\alpha}')$, which can
be easier described.

Free loop spaces are related with our problem because of the following well known result.
See also the work of Lin [L] for the following considerations.

\vspace{0.2cm}

\noindent \textbf{Proposition.} \ \textit{For each map $f: P=\coprod_rS^1\rightarrow M$
there is a natural isomorphism of groups:
$$\pi_1(\tilde{\mathcal{M}},f)\rightarrow \times_{i=1}^r\pi_1(LM,f_i),$$
where $f_i$ is the restriction of $f$ to the $i$-th component.}

\vspace{0.2cm}

By 4.7 we can restrict to sets of generators. So the submodule
$W''(r)$ can be computed from the contributions of the groups
$\pi_1(LM,f_i)$ for components $f_i$ of standard links $K_{\alpha
}$. We have to compute the image $W''(r)$ in
$RC(r)/(W(r-1)+W'(r))$ or the corresponding lift to $RC(r)$. Each
element of $\pi_1(LM,f_i)$ determines a homotopy class of loops in
$K_{\alpha }$ (after projecting to $\mathcal{M}$), defined by a
self homotopy of a component  with homotopy class $f_i$. We
choose a representative path $\gamma \in \mathcal{P}_{\alpha }$
and compute $s_{if}(\gamma )$. The result corresponds to the
definition of $\chi $ by the $\Theta$-construction. Because of
condition 3 for geometric standard links, the contributions for
different components but with the same homotopy class coincide.

\vspace{0.5cm}

\centerline{\textsc{6.\ EXPANSION OF RELATIONS}}

\vspace{0.3cm}

We need a workable algebraic description of the fundamental group of free loop spaces. The results for $\pi_2(M)=0$ are in [L].

\vspace{0.2cm}

\noindent \textbf{Lemma 6.1.} \ {Let $f: S^1\rightarrow M$ be a map with homotopy class $[f]\in \pi_1(M)$. Then there
is the exact sequence
of groups $$\CD
\pi_2(M)@>{\pi_f^*}>>\pi_1(LM,f)@>{i_f^*}>> Z([f])\rightarrow 1,
\endCD$$ where for each element $b\in \pi_1(M)$ we let $Z(b)$ denote the
centralizer of $b$. The basepoint for the definition of $\pi_2(M)$ and $Z(b)\subset \pi_1(M)$ is
$f(1)$.}

\vspace{0.2cm}

\noindent \textsl{Proof.} \ The homomorphism $\pi_f^*$ is defined
as follows: Consider the map $\bar{f} :S^1\vee S^1\rightarrow M$,
which restricts to $f$ on $S^1\times \ast$ and is constant on
$\ast \times S^1$. Consider the attaching map of the $2$-cell $c:
S^1\rightarrow S^1\vee S^1$. Then $\bar{f}\circ c$ is canonically
null-homotopic in $M$ (contracting the path determined by $f$ to
the initial point). So each element of $\pi_2(M)$ determines a
unique up to homotopy extension of $\bar f$ over the attached
$2$-cell. Thus we have defined a map $S^1\times S^1\rightarrow
M$, which represents an element of $\pi_1(LM,f)$. It is easy to
see that this actually defines a homomorphism with image the
kernel of $i_f^*$. Here $i_f^*$ is defined by composition of the
natural map $\pi_1(LM,f)\rightarrow [S^1\times S^1,M]$ with the
map induced by the inclusion $\ast \times S^1\subset S^1\times
S^1$. Here $[S^1\times S^1,M]$ denotes the set of all based
homotopy classes of maps $S^1\times S^1\rightarrow M$. The
exactness follows by obstruction theory. $\square$

\vspace{0.2cm}

For $\alpha \in \mathfrak{b}(M)$ and $f$ representing a component of
$K_{\alpha }$ with homotopy class $a$, define
$$\iota_{f,\alpha}:\pi_1(LM,f)\rightarrow \mathbb{Z}$$
by $\iota_{f,\alpha}(h)=\iota_f(h,\alpha\setminus a)$
(compare section 1). Let $f_1,\ldots ,f_r$ be maps representing
the components of $K_{\alpha }$. For $h\in
\times_{i=1}^r\pi_1(LM,f_i)$ define
$$\varepsilon (h)=\sum_{i=1}^r\iota_{f_i,\alpha}(h_i).$$

\vspace{0.2cm}

The proposition in section 5 identifies $\varepsilon $ with the
the index function on $\pi_1(\tilde{\mathcal{M}},K_{\alpha})$
defined in 4.1.

\vspace{0.2cm}

A homomorphism $\varepsilon :G\rightarrow \mathbb{Z}$ is called an \textit{indexed group}. Let $A$ be a $\mathbb{Z}[q^{\pm 1}]$-module and $\varepsilon :G\rightarrow
\mathbb{Z}$ be an indexed group. A map $\chi : G\rightarrow A$ is  called a
\textit{homomorphism} if $\chi (1)=0$ and  $\chi (g_1g_2)=\chi
(g_1)+q^{2\varepsilon(g_1)}\chi (g_2)$.
For an indexed group $\varepsilon :G\rightarrow \mathbb{Z}$ and $i$ a non-nengative integer let $G^i:=\varepsilon^{-1}(i)$.

\vspace{0.2cm}

\noindent \textbf{Lemma 6.2.} \textit{Let $\varepsilon : G\rightarrow \mathbb{Z}$ be an indexed group and $\chi : G\rightarrow A$ be a homomorphism. Then $\chi $
restricts to a group homomorphism
$\chi^0: G^0\rightarrow A$, which factors
factors through $G^0/[G^0,G^0]$. Here $[G^0,G^0]$ is the commutator subgroup.}

\vspace{0.2cm}

\noindent \textbf{Lemma 6.3.} \ \textit{For each indexed group $G$
and exact sequence of groups:
$$\CD H@>{j}>>G@>>>Q@>>>1 \endCD$$
there exists the commutative diagram of groups:}
$$\CD
H^0 @>{j|}>> G^0 @>>> Q^0 @>>>1 \\
@VV{\subset}V @VV{\subset}V @VV{\subset}V \\
H @>{j}>> G @>>> Q @>>> 1 \\
 @. @VV{\varepsilon}V @VV{\bar{\varepsilon}}V @. \\
 @. \mathbb{Z} @>>> \mathbb{Z}/(\varepsilon \circ j(H)) @.
\endCD
$$ \textit{with the canonical isomorphism
$G^0/(j|(H^0))\rightarrow Q^0$. Here $H^0:=ker(\varepsilon \circ
j)$, $\bar{\varepsilon }$ is the canonical homomorphism induced by
$\varepsilon $ and $Q^0=ker(\bar{\varepsilon })$. Moreover, there
is the exact sequence $$\CD H^0@>>>G^0/[G^0,G^0]@>>>Q^0/[Q^0,Q^0]
@.  @>>>0\endCD$$}

\noindent \textsl{Proof.} \ This is proved by standard diagram chase arguments.
$\square$

\vspace{0.2cm}

We deduce the following consequence:

\vspace{0.2cm}

\noindent \textbf{Proposition 6.4.}\ \textit{Let $a\in \pi_1(M)$ be contained in $\alpha $. Then the module generated by $\chi (T(a))$ is also generated by the images of $\chi $ on generators of abelian groups $$\pi_2(M), \ Z(a)/Z(a)^0 \ \  \text{and} \ \
Z(a)^0/[Z(a)^0,Z(a)^0].$$}
\noindent \textsl{Proof.} \ Apply 6.3 and 6.2 to the exact sequence of 6.1. Define $\chi $ on $Z(a)$ by a section of maps
\textit{maps} $Z(a)\rightarrow \pi_1(LM,f)$, which maps into
$\pi_1(LM,f)^0$. $\square$

\vspace{0.2cm}

Now we describe first order expansions of the elements
$\Delta (\alpha ,(a,b),g)$ and $\Theta (\alpha ,a,h)$ from section 1 and  compute important examples. The ordering and basing of $K_{\alpha}$
lifts $\alpha =\langle \alpha_1, \ldots ,\alpha
_r\rangle$ to $(a_1,\ldots ,a_r)$.
Then for each $g\in \pi_1(M)$ define
$$\alpha (i,j;g):=\langle (a_iga_jg^{-1})^{\circ} ,\alpha_1,\ldots ,\hat{\alpha_i},\ldots,\hat{\alpha_j},
\ldots,\alpha_r\rangle.$$

\vspace{0.2cm}

Assume that the based first based component of $K_{\alpha }$ has homotopy class $a\in \pi_1(M)$. For $2\leq i\leq r$ and $h\in T(a)$ let $\varepsilon_i=\varepsilon_i (h)$ denote the intersection number of the map of the torus given by $h$ with the $i$-th component of $K_{\alpha }$. Note that $\varepsilon (h)=\sum_{i=2}^r\varepsilon_i (h)$.

\vspace{0.2cm}

\noindent \textbf{Proposition 6.5.} \ \textit{Let $|\alpha |=r$ and $h\in
T(a)$ for some $a$ in $\alpha $ as above. Then we can choose the expansion
in the $\Theta$-construction so that
$$\Theta (\alpha, a,h)=(q^{2\varepsilon (h)}-1)\alpha +zp_1(h).$$ Moreover,
for $2\leq i\leq r$ there exist natural numbers $n_i\geq 0$ and for $0\leq j\leq
n_i$ numbers $\varepsilon_{i,j}\in \{+1,-1\}$
with
$\sum_{j=0}^{n_i}\varepsilon_{i,j}=\varepsilon_i(h)$. Also there exist elements
$g_{i,j}\in
\pi_1(M)$ and integer numbers $\eta_{i,j}$ (all numbers and group elements depending on $h$)
such that $$\CD p_1(h)&=\sum_{i=2}^r q^{2(\varepsilon_2+\ldots
+\varepsilon_{i-1})}\sum_{j=0}^{n_i} \varepsilon_{i,j}q^{2(\varepsilon_{i,1}+\ldots
+\varepsilon_{i,j-1})+\varepsilon_{i,j}+\eta_{i,j}}\alpha(1,i;g_{i,j})\\
&+z\tilde{p}_1(h)\endCD$$ for some $\tilde{p}_1(h)\in SR\pi(M)$.}

\vspace{0.2cm}

\noindent \textsl{Proof} The self-homotopy of $f$ with $[f]=a$
defines $f': S^1\times S^1\rightarrow M$. Assume that $f'$
intersects $K_{\alpha ,a}$ transversely. Recall that
$f'|S^1\times *=f$. Homotope the inclusion $S^1\times *\subset
S^1\times S^1$ by pushing arcs successively across intersection
points. First move arcs across intersections with the second
component of $K_{\alpha }$, then with the third component, and so
on (the first component is deformed). Finally we have deformed
$S^1\times *$ across all intersection points and can isotope the
result back to the original inclusion without further
intersections. The composition of this deformation with $f'$
defines a self-homotopy of $f$ homotopic to the given one.
$\square$

\vspace{0.2cm}

\noindent \textbf{Example 6.6.} \ Let
$M=S^1\times S^1\times S^1$ with $\pi_1(M)\cong \mathbb{Z}^3$
generated by elements $b_i$, $i=1,2,3$. Represent $M$ as quotient
of $I^3$ by identification of opposite sides. Let $\alpha
=\langle b_1,b_2,b_3 \rangle$. The standard link $K_{\alpha }$
can be represented in $I^3$ by intervals parallel to the
faces of $I^3$. Consider the self homotopy of $K_{\alpha }$, which is
determined by $h\in \pi_1(LM,f_1)$ with $[f_1]=b_1$
as follows: Keep fixed all but the first component with homotopy
class $b_1$. The first component crosses the $b_2$-component and
then the link is isotoped back to the original standard link in
the obvious way. Note that the resulting map of the torus
restricted to $S^1\vee S^1$ has homotopy class given by $b_1$ and
$b_3$. We have $\varepsilon (h)=\varepsilon_2(h)=1$ and compute:
$$\Theta(\alpha,a,h)= (q^2-1)\alpha +zq\langle b_1b_2, b_3 \rangle.$$
Now, after the deformation above, let the first component
additionally crosses the third component with negative sign. Let
$h'$ be the resulting element of $\pi_1(LM,f_1)$. The restriction
of the corresponding map of the torus to $S^1\vee S^1$ has
homotopy classes given by $b_1$ and $b_3b_2^{-1}$. We have
$\varepsilon (h')=0$ and $\varepsilon_2(h')=1$,
$\varepsilon_3(h')=-1$, and compute:
$$\Theta(\alpha,a,h')=zq(\langle b_1b_2, b_3\rangle - \langle b_1b_3, b_2 \rangle).$$
The element $\sigma (\langle b_1b_2, b_3 \rangle ) - \sigma (\langle
b_1b_3, b_2 \rangle)$ is a torsion element in
$\mathcal{H}(S^1\times S^1\times S^1)$ by 2.4.

\vspace{0.2cm}

\noindent \textbf{$\Theta$-Construction for the image of $\pi_f^*$:} \ Let $f:S^1\rightarrow M$ represent a component of $K_{\alpha }$ and $[f]=a$.
By the homomorphism
$\pi_f^*: \pi_2(M)\rightarrow \pi_1(LM,f)$ each element $h\in \pi_2(M)$
defines a self homotopy
of $f$. Then the construction of $\Theta (\alpha, a,\pi_f^*(h))$
for $h\in \pi_2(M)$ can be described as follows: Let $f': S^2\rightarrow M$
represent $h$ with $f'(*)=f(1)$. We can change $f'$ by homotopy such that $f'$ restricts to
$f$ on some closed interval $I$ containing $*$ in its interior by using some homeomeorphism $j$
between $I$ and the closed half circle $I'$ containing $1\in S^1$. Now let $H$ be a homotopy of
the interval $I$ across the $2$-sphere fixing $\partial I$ such that the image of the homotopy
covers the $2$-sphere with degree 1. This defines a self homotopy of $f$ by $f|I'= f'\circ H
\circ j^{-1}$ and $f$ is constant on $S^1\setminus I'$.

\vspace{0.2cm}

Recall the action of $u$ on $\mathfrak{b}(M)$ defined in the proof of 2.3.

\vspace{0.2cm}

\noindent \textbf{Proposition 6.7.} \ \textit{For $\alpha \in \mathfrak{b}(M)$ and
$h\in \pi_2(M)$, the $\Theta$-construction can be chosen such that
$$\Theta (\alpha ,a,\pi_f^*(h))= (u(q^{2\varepsilon(h)}-1)+zq \frac{1-q^{2\varepsilon (h)}}{1-q}) \alpha
+z^2p_0(h)$$ for some $p_0(h)\in \mathcal{R}\hat{\mathfrak{b}}(M)$.}

\vspace{0.2cm}

\noindent \textbf{Example 6.8.} \ Let $M=S^2\times S^1$ and let $h\in \pi_2(M)$
be a generator. We represent $h$ by the standard embedded sphere $S:=S^2\times \{1\}
\subset S^2\times S^1$.
Let $a\in \pi_1(M)\cong \mathbb{Z}$ be a generator und $\alpha =\langle a,a,a^{-1}\rangle$. We
assume that the ordering and basing of $K_{\alpha }$ lifts $\alpha $ to $(a,a,a^{-1})$. We can assume that the first component
intersects $S$ in a closed interval in the complement of the two transverse intersection points
of the last two components with $S$. Then we apply the $\Theta$-construction for a
self-homotopy of the first component. Note that $\varepsilon (h)=0$. We compute $\chi (h)$ in this case. By expanding at the crossing points we find that the
contributions add up to $zq(K-K')$ for two links with $\omega (K)=\langle
a^2,a^{-1}\rangle $ and $\omega (K')=\langle 1,a \rangle $. If we expand $K,K'$ in terms of
standard links we compute
$$\Theta(\alpha,a,\pi_f^*(h))=zq(q^{\eta_1} \langle a^2,a^{-1}\rangle - q^{\eta_2}
\langle 1,a \rangle)+z^2\tilde{p}_1(h)$$
for numbers $\eta_1,\eta_2\in \mathbb{Z}$. It follows that $S^2\times S^1$ has torsion.

\vspace{0.3cm}

Now we consider the $\Delta$-construction.

\vspace{0.2cm}

\noindent \textbf{Proposition 6.9.} \ \textit{Let $\alpha \in \mathfrak{b}(M)$ and $(a,b)$ a pair of elements (w.\ r.\ lifting the first two components of the ordered link) in $\alpha $ and $g\in \pi_1(M)$. Then
there exists $\delta (g)\in \mathbb{Z}$ and $p_2(g)\in
SR\hat{\pi }(M)$ (containing only elements of $\mathfrak{b}(M)$ of
length $\leq |\alpha |-2$) such that
$$\Delta (\alpha ,(a,b),g)
=z(q^{-1}-q)(q^{2\delta (g)}\alpha (1,2;g)-\alpha
(1,2;1))+z^2(q^{-1}-q)p_2(g).$$}

\vspace{-0.3cm}

\noindent \textbf{Example.}
Let $M=S^1\times S^1\times S^1$ and let
$b_i$, $i=1,2,3$ be the natural generators of
$\pi_1(M)$. Let $\alpha =\langle b_1,b_1^{-1},b_2
\rangle$ and $c=(b_1,b_1^{-1})$. Then $\alpha
(1,2;1)=\alpha (1,2;b_3)=\langle 1,b_2 \rangle$ because $\pi_1(M)$
is commutative. We can choose the basings such that
$K_{\alpha }^c$ is the standard link
with wrapping invariant $\langle 1,b_2 \rangle$. The result of
attaching the band determined by $b_3$ is a link, which can be
deformed into the standard link through a single crossing change
with the component of homotopy class $b_2$. We compute
$$\tilde \Delta (\alpha ,c,
b_3)= (1-q^2)\langle 1,b_2 \rangle -qz\langle b_2\rangle.$$ The
element $(1-q^2)\langle 1,b_2\rangle -qz\langle b_2\rangle $ is
not trivial in $\mathcal{H}(M)$. In order to show this apply the
homomorphism $z\mapsto 0$ and consider the resulting element of
$\mathcal{L}(M)$. The result follows because $\lambda \langle
1,b_2\rangle =0$ (compare theorem 1.2). In fact
$\iota_f(b,b_2)=0$ for $b\in Z(1)\cong \pi_2(M)=0$ and $f$
representing the trivial homotopy class, and $\iota_{f'}(b,1)=0$
for all $b\in Z(b_2)$ and $f'$ representing $b_2$ (since
intersections with a trivial homology class vanish). Thus
$\mathcal{H}(S^1 \times S^1\times S^1)$ has torsion coming from
the $\Delta$-construction, even though the fundamental group is
abelian.

\vspace{0.5cm}

\centerline{\textsc{7.\ REDUCTION OF STRUCTURE SETS}}

\vspace{0.3cm}

We prove that $\chi $ is trivial in certain situations. This
appears in connection with \textit{cyclic} and \textit{peripheral} elements for both
the $\Delta$- and $\Theta$-construction.

\vspace{0.2cm}

Throughout let $T_1,\ldots T_n$ be the tori in
$\partial M$. We choose paths from $*\in M$ to fixed points in $T_I$ for each $i$,
thus define homomorphisms $j_* : \pi_1(T_i)\rightarrow \pi_1(M)$ for $1\leq i\leq
n$.

\vspace{0.2cm}

\noindent \textbf{Proposition 7.1.} \ \textit{Let $M$ be a $3$-manifold with $\pi_1(M)$ cyclic. Then for all
$g\in D(a,b)$ and $(a,b)$ a pair in $\alpha $ we have $\Delta (\alpha ,(a,b),g)=0$ .}

\vspace{0.2cm}

\noindent \textsl{Proof.} Let $T\hookrightarrow int(M)$ be an
embedded solid torus, which induces an epimorphism of fundamental
groups. Let $K_{\alpha }$ be a link and let $(a,b)$ be the
homotopy classes of the first two components. Consider a band
determined by $g$ (compare the $\Delta$-construction). We only
need to consider the arc (joining the two components) up to
\textit{homotopy} with endpoints on the components. So we can
assume that the first two components are contained in $T$ while
all the other components are in $M\setminus int(T)$. There could
be intersections of different components during such a homotopy.
But we can assume that the smoothings will contribute only terms,
which have already been taken care of by induction, similar to
the proof of 2.1 in section 5. The band will be homotoped along.
Then we can homotope the band into $T$, again intersections with
the link can be neglected by the arguments from section 6. But
since $\pi_1(T)$ is abelian the results of the standard banding
and the one determined by $g$ are homotopic in $T$. So the skein
difference actually vanishes. In general the resulting link is
not a standard link. But for both ways of smoothing it is the same
link, so the expansions in terms of standard links can be assumed
equal. $\square$

\vspace{0.2cm}

\noindent \textbf{Proposition 7.2.} \textit{Let $M=S^1\times
S^1\times I$. Then $\Delta (\alpha, a,g)=0$ for all $g\in D(a,b)$
and for all pairs $(a,b)$ in $\alpha $. }

\vspace{0.2cm}

\noindent \textsl{Proof.} \ The idea is similar as in 8.1. We can
assume that the two distinguished components and the attaching
band are contained in $S^1\times S^1\times [0,\frac{1}{2}]$ while
the remaining components of $K_{\alpha }$ are contained in
$S^1\times S^1\times [\frac{1}{2},1]$. Because $\pi_1(S^1\times
S^1)$ is abelian, the links which result from attaching the bands
are homotopic to each other in $S^1\times S^1\times
[0,\frac{1}{2}]$. The rest of the argument is similar to 7.1.
$\square $

\vspace{0.2cm}

The two geometric observations above generalize to arbitrary structure sets.
\begin{enumerate}
\item Assume that $a,b \in \langle h \rangle$ for some $h\in \pi_1(M)$.
Then let $R_1(a,b)$ denote the set of cosets of elements $h^i$
for all $i\in \mathbb{Z}$.
\item Assume that $a$ or $b$ is conjugate to the image of a torus
$T_i\subset \partial M$, e.\ g.\
$h^{-1}ah\in j_*(\pi_1(T_i))$ for some $h\in \pi_1(M)$ . Let
$\Pi_i$ be the union of all corresponding subgroups $h^{-1}j_*(\pi_1(T_i))h$.
Then let
$R_2(a,b)$ be  the set of cosets $\langle a\rangle g \langle b\rangle$ with $g\in \Pi_i$
for some $1\leq i\leq n$.
\end{enumerate}

\noindent \textbf{Theorem 7.3} \ \textit{
$\chi_{\alpha}(R_1(a,b)\cup R_2(a,b))=0$ for all $\alpha $ and $(a,b)$ in $\alpha $.}

\vspace{0.2cm}

\noindent \textsl{Proof.} \ For cosets from $R_1(a,b)$ we use the torus argument
as before. We can homotope both $a,b$ into the torus determined by $h$ as before. If the
band is determined by some element in $\pi_1(M)$, which can be homotoped into the torus,
then the argument follows 7.1. For $R_2(a,b)$ assume that $a$ is freely homotopic into some boundary torus $T_i$ and let $h$ be the conjugating element. The band is determined by some element in $\pi_1(M)$ and by assumption can be homotoped, such that it runs inside a collar neighbourhood of the torus and otherwise follows the standard band. The rest of the argument is like in 7.2. (use the basepoint in $T_i$). $\square$.

\vspace{0.2cm}

Next we consider $T(a)$ for some $a\in \pi_1(M)$. Again we can
reduce in peripheral and a cyclic situations.

\vspace{0.2cm}

\noindent \textbf{Proposition 7.4.} \ \textit{Let $M=F\times I$ for $F$ a compact
oriented surface. Then $\chi _{\alpha }(T(a))=0 $  for all $a$ in $\alpha $.}

\vspace{0.2cm}

\noindent \textsl{Proof.} We can homotope singular tori into
$F\times [0,\frac{1}{2}]$ while keeping the rest of the link in
$F\times [\frac{1}{2},1]$. In general such a homotopy of a
singular torus is not a homotopy of loops in $LM$ rel $f$. Let
$\gamma $ denote a loop in $\mathcal{M}$ in $K_{\alpha }$ in
almost general position, which represents the given element in
$\pi_1(LM,f)$. Consider the homotopy as a map $S^1\times
I\rightarrow \mathcal{M}$, which is in almost general position on
the boundary. By 3.1 we can approximate by a map $S^1\times
I\rightarrow \mathcal{M}_2$ in almost general position. Moreover
we can assume that the restriction to $\{1\}\times I$ is in
almost general position. Let $\gamma_0$ be the path defined by
the restriction to $\{1\}\times I$ and let $\gamma '$ be the loop
in $\gamma_0(1)$, which is determined by the restriction to
$S^1\times \{1\}$. We have to compare the two terms
$s_{if}(\gamma )$ and $s_{if}(\gamma_0\gamma' \gamma_0^{-1})$.
Since there are no singular parameters on $\gamma '$ we know that
$s(\gamma ')=0$, in particular $\varepsilon (\gamma ')=0$. It
follows that also $\varepsilon (\gamma_0\gamma
'\gamma_0^{-1})=\varepsilon (\gamma )=0$. Because of 4.6 it
suffices to consider the difference between $s(\gamma )$ and
$s(\gamma_0\gamma'\gamma_0)$, which has been computed in 4.5. It
follows that $s_{if}(\gamma )=0$ and so the assertion. $\square $

\vspace{0.2cm}

\noindent \textbf{Proposition 7.5.} \ \textit{Assume that
$\pi_2(M)=0$. Let $a\in \pi_1(M)$, $a=b^k$ for some $k\in
\mathbb{Z}$ and and $h\in Z(a)\cap \langle b \rangle$. Then $\chi
(h)=0$ (use the isomorphism $Z(a)\cong \pi_1(LM,f)$ for
$[f]=a$ from 6.1).}

\vspace{0.2cm}

\noindent \textsl{Proof.} \ Find a knot in $M$ with homotopy
class $b$. There is a regular neighbourhood $T$ disjoint from the
link $K_{\alpha}$. Because of the assumptions we can first
homotope the map of the torus restricted to its $1$-skeleton $T$.
But then because of $\pi_2(M)=0$ the map of torus is homotopic
into $T$, because $T$ is disjoint from other components of
$K_{\alpha }$.  Since $\pi_2(M)=0$ the map of torus can be
homotoped into $T$. Then we can approximate the resulting map in
$\mathcal{M}$ by a map in almost general position in $T$. The
arguments from 7.4 also apply in this case. The intersections
with the components of $K_{\alpha ,a}$ are trivial so $\chi
(h)=0$. $\square $

\vspace{0.2cm}

Again the two ways of reductions generalize.

\begin{enumerate}
\item Note that $\langle a \rangle \subset Z^0(a)$. If $a=b^k$ for some $k\in \mathbb{Z}$ then
let $\langle a \rangle _c$ denote the subgroup
$Z(a)\cap \langle b \rangle$.
\item Assume that $a$ is conjugate to some element in $j_*(\pi_1(T_i))$ for
some $1\leq i\leq n$.
Then let $g\in \pi_1(M)$ such that $gag^{-1}\in j_*(\pi_1(T_i))$ for
some $i$, $1\leq i\leq n$. Let $\Pi_i:=gj_*(\pi_1(T_i))g^{-1}\subset
Z^0(a)$, and let $\Pi $ be the product of the $\Pi_i$.
\end{enumerate}

\noindent \textbf{Theorem 7.6.} \ \textit{The image of the
homomorphism $\chi $ is generated by its images on generating sets
of the three abelian groups: $$\pi_2(M), \ Z^0(a)/( \langle a
\rangle _c [Z^0(a),Z^0(a)] \Pi ) \ \text{and} \ Z(a)/Z^0(a).$$}

\vspace{-0.3cm}

\noindent \textsl{Proof.} By 6.3 it suffices to consider sets of
generators of $\pi_2(M)$, $Z(a)/Z^0(a)$ and $Z^0(a)/[Z^0(a),Z^0(a)]$.
By the argument from 7.5, $\chi $ is trivial on
$\langle a \rangle_c$. A singular torus, which is determined by
some element of $\Pi $, can be homotoped into a collar
neighborhood of a boundary torus. The rest of  the argument is
analogous to the proof of 7.4. $\square $.

\vspace{0.5cm}

\centerline{\textsc{8.\ PROOFS OF THEOREM 1.1 AND THEOREM 1.2}}

\vspace{0.3cm}

For standard results about $3$-manifolds we refer to [He].

\vspace{0.2cm}

Theorem 1.2 follows from 1.3 and the results of section 6.
For the proof of 1.1 b), by the
universal coefficient result [P-3] and right-exactness of the tensor product,
it suffices to show $\chi (T(a))=0$,
if $M$ is atoroidal and $\pi_2(M)=0$. By [He] 9.13,
$f_*(\pi_1(S^1\times S^1))$ is cyclic or $\mathbb{Z}^2$ for each
map $f$ of a torus in $M$. If it is cyclic then apply 7.5 or 7.6.
If the torus map induces an injection of fundamental groups
then by assumption the corresponding element $h\in Z(a)$ is contained
in a subgroup $\Pi_i$, compare 7.6. So $\chi (T(a)=0$
for all $a$ in $\alpha $. So we only have to
prove theorem 1 a).

\vspace{0.2cm}

For each $3$-manifold $M$ let $\mathcal{P}(M)$ denote the \textit{Poincare
associate} of $M$, which is defined by first capping off all
$2$-spheres in the boundary of $M$ and then replacing a possibly
fake $3$-cell in $M$ by a standard $3$-cell (w.\ r.\ we can
assume that there is at most one fake $3$-cell in $M$).

\vspace{0.2cm}

\noindent \textbf{Lemma.} \ \textit{For each $3$-manifold $M$ there
is a natural isomorphism $$\mathcal{H}(M)\cong \mathcal{H}(\mathcal{P}(M)),$$
induced by inclusions.}

\vspace{0.2cm}

\noindent \textsl{Proof.} The inclusions $M\hookrightarrow \hat
M$ and $M\setminus e\hookrightarrow \mathcal{P}(M)$ ($e$ is the
interior of a fake $3$-cell in $\hat{M}$) induce isomorphisms
because of the same reason. We consider the inclusion
$\hat{M}\setminus e\hookrightarrow \hat{M}$. Since the induced
map of fundamental groups is an isomorphism it also induces
bijections between the sets of conjugacy classes. So we can
choose a geometric set of standard links in $\hat{M}\setminus e$,
which is also a set of geometric standard links for $\hat{M}$
(see 2.2). But band constructions and homotopies of singular tori
can be assumed in $\hat{M}\setminus e$. So the homomorphisms
$\chi $ for $\hat{M}\setminus e$ and $\hat{M}$ are compatible.
The result follows by comparison of the presentation sequences in
2.2. $\square $.

\vspace{0.2cm}

Assume that $\pi_1(M)$ not abelian. Let $a,g \in \pi_1(M)$ be
elements with \newline $aga^{-1}g^{-1}\neq 1$ and $\alpha :=\langle
a^{\circ},(a^{-1})^{\circ} \rangle $. We can choose basings of
the standard links such that $a^{\circ}$ lifts to $a$. Then
$\alpha (1,2;1)=1^{\circ}$ and $\alpha
(1,2;g)=(aga^{-1}g^{-1})^{\circ}\neq 1^{\circ}$. It follows from
1.3 and 2.4 that $\sigma $ maps $\tilde{\Delta }(\alpha, (a,a^{-1}), g)$ to
a torsion element of $\mathcal{H}(M)$.

\vspace{0.2cm}

Assume that $\pi_1(M)$ is abelian. Then by [He, 9.13],  $\pi_1(M)$
is cyclic or free abelian of rank 2 or 3. Moreover, if $M$ is
free abelian of rank 2 resp.\ 3 then $\mathcal{P}(M)$ is is
homeomorphic to $S^1\times S^1\times I$ or $S^1\times S^1\times
S^1$ (see [He 5.2 and 11.11]). If $M$ is free abelian of rank 1
then $\mathcal{P}(M)$ is homeomorphic to $S^1\times D^2$ or
$S^1\times S^2$. Note that the condition $2b_1(M)=b_1(\partial
M)$ ($\Leftrightarrow$ $M$ is submanifold of a rational homology
$3$-sphere [K-2]) is unchanged if fake $3$-cells are replaced by
standard $3$-cells. In particular, among the $3$-manifolds with
$\pi_1(M)$ abelian precisely those with finite cyclic fundamental
group or $\mathcal{P}(M)$ homeomorphic to $S^1\times D^2$ or
$S^1\times S^1\times I$ are those, which are submanifolds of
rational homology $3$-spheres. Note that by the lemma above we
can replace $\hat{M}$ by $\mathcal{P}(M)$. So assume that
$\pi_1(M)$ is cyclic, but $\mathcal{P}(M)$ is not homeomorphic to
$S^1\times S^2$. Then it follows from 7.1, 7.4 and 7.5 that $\chi
$ is the trivial homomorphism and $\mathcal{H}(M)$ is free. If
$\mathcal{P}(M)$ is homeomorphic to $S^1\times S^1\times I$ then
the same follows from 7.2 and 7.4. In 6.6 and 6.8 we have shown
that $\mathcal{H}(S^1\times S^1\times S^1)$ and
$\mathcal{H}(S^1\times S^2)$ have torsion. This completes the
proof of 1.1. $\square$

\vspace{0.2cm}

\noindent \textbf{Remark.} \ It is immediate from 6.5 that
$\mathcal{C}(M)$ has torsion, if $M$ contains a singular torus,
which has non-trivial intersection number with some oriented loop
in $M$. But in general it is difficult to construct torsion in
$\mathcal{C}(M)$ [K-4].

\vspace{1cm}

\centerline{\textsc{REFERENCES}}

\begin{itemize}
\item[[B]]
\textsc{D.\ Bullock}. \textsl{Rings of
$Sl_2(\mathbb{C})$-characters and the Kauffman bracket Skein
module}. Comentarii Math.\ Helvetici 72, 1997, 619--632
\item[[B-F-K]]
\textsc{D.\ Bullock, C.\ Frohman, J.\ Kania-Bartoszynska}.
\textsl{Understanding the Kauffman bracket skein module}.
Journal of Knot Theory and Its Ramifications, Vol 8, No. 3, 1999, 265--277
\item[[H]]
\textsc{N.\ Habegger}. \textsl{Link homotopy in simply connected
$3$-manifolds}. \newline AMS/IP Studies in Mathematics,
Proceedings of the Georgia International Topology Conference
1997, 118--123
\item[[He]]
\textsc{J.\ Hempel}.
\textsl{$3$-manifolds}, Annals of Mathematical Studies 86.
Princeton University Press 1976
\item[[H-P]]
\textsc{J.\ Hoste, J.\ Przytycki}.
\textsl{Homotopy skein modules of orientable $3$-manifolds}.
Math. Proc.\ Cambr.\ Phil.\ Society 108, 1990, 475--488
\item[[J-S]]
\textsc{W.\ H.\ Jaco, P.\ B.\ Shalen}.
\textsl{Seifert fibred spaces in $3$-manifolds}.
Memoirs of the AMS 220, Vol. 21, 1979
\item[[K-1]]
\textsc{U.\ Kaiser}.
\textsl{Link homotopy and skein modules of $3$-manifolds}.
Geometric Topology, ed. C.\ Gordon, Y.\ Moriah, B.\ Wajnryb, Contemp.\ Math.\ 164, 1994,
59--77
\item[[K-2]]
\textsc{U.\ Kaiser}.
\textsl{Link theory in manifolds}.
Lecture Notes in Math.\ 1669, Springer Verlag, 1997
\item[[K-3]]
\textsc{U.\ Kaiser}
\textsl{Homology boundary links and fusion constructions}.
Osaka J.\ of Math.\  29, 1992, 573--593
\item[[K-4]]
\textsc{U.\ Kaiser}.
\textit{Torsion in homotopy skein modules}.
in preparation
\item[[Ka-1]]
\textsc{E.\ Kalfagianni}
\textsl{Finite type invariants for knots in $3$-manifolds}.
Topology 37, no.3,  1998, 539--574
\item[[Ka-2]]
\textsc{E.\ Kalfagianni}.
\textsl{Power series link invariants and the Thurston norm}.
Topology and its Applications 101, 2000, 107--119
\item[[Ka-L]]
\textsc{E.\ Kalfagianni, X.\ S.\ Lin}.
\textsl{The HOMFLY polynomial for links in rational homology spheres}.
Topology 38, no. 1, 1999, 95--115
\item[[KiLi]]
\textsc{P.\ Kirk, C.\ Livingston}.
\textsl{Knot invariants in $3$-manifolds, essential tori and the first
cohomology of the free loop space on $M$}.
preprint 1997
\item[[L]]
\textsc{X.\ S.\ Lin}.
\textsl{Finite type link invariants of
$3$-manifolds}.
Topology 33, no.1, 1994, 45--71
\item[[P-1]]
\textsc{J.\ Przytycki}.
\textsl{Skein modules of $3$-manifolds}.
Bull.\ Ac.\ Pol.\ Math.\ 39(1-2), 1991, 91-100
\item[[P-2]]
\textsc{J.\ Przytycki}.
\textsl{$q$-analog of the homotopy skein module of links}.
preprint, Knoxville 1991
\noindent [P-3] J.\ Przytycki, \textsl{Homotopy and $q$-homotopy
skein modules of $3$-manifolds: an example in Analysis Situs},
preprint 1999
\item[[P-4]]
\textsc{J.\ Przytycki}.
\textsl{Skein module of links in a handlebody}.
Topology 90, Proc.\ of the Research Semester in Low
Dimensional Topology at OSU, Editors: B.\ Apanasov, W.\ D.\
Neumann, A.\ W.\ Reid, L.\ Siebenmann, DeGruyter Verlag, 1992,
315--342
\item[[P-5]]
\textsc{J.\ Przytycki}.
\textsl{A q-analogue of the first homology group of a $3$-manifold}.
Contemporary Mathematics Vol.\ 214, 1998, 135--144
\item[[P-S]]
\textsc{J.\ Przytycki, A.\ Sikora}.
\textsl{On skein algebras and $Sl_2(\mathbb{C})$-character varieties}.
Topology 39, 2000, 115--148
\item[[T]]
\textsc{V.\ Turaev}.
\textsl{Skein quantization of Poisson algebras of loops on surfaces}.
Ann.\ Scient.\ Ec.\ Norm.\ Sup.\  4 (24), 1991, 635--704
\end{itemize}

\end{document}